\documentclass[11pt,titlepage]{article}

\usepackage[all]{xy}
\usepackage{amssymb,latexsym,amsfonts,amsmath, mathrsfs}
\usepackage{graphicx,diagrams}
\usepackage{comment}
\usepackage{epsfig}
\usepackage{color}
\usepackage[scale=0.8, centering]{geometry}

\newtheorem{theorem}{Theorem}[section]

\newtheorem{proposition}[theorem]{Proposition}
\newtheorem{corollary}[theorem]{Corollary}

\newtheorem{definition}[theorem]{Definition}
\newtheorem{example}[theorem]{Example}

\newtheorem{remark}[theorem]{Remark}
\numberwithin{equation}{section}

\renewcommand{\Re}{{\mathbb{R}}}
\newcommand{\R}{{\mathbb{R}}}

\newcommand{\N}{{\mathbb{N}}}

\newcommand{\e}{\mathbf{e}}

\begin{document}
%
% paper title
\title{
\Huge{To sample or not to sample: \\Self-triggered control for nonlinear systems}\\[2cm]
%\author{Adolfo Anta}\\[2cm]% <-this % stops a space
\Large{Adolfo Anta and Paulo Tabuada} \\[8cm]
\begin{flushleft}
\small{This research was partially supported by the National Science Foundation EHS award 0712502 and a scholarship from Mutua Madrile\~na Automovilista.}\\[0.5cm]
\small{A. Anta and P. Tabuada are with the Dept. of Electrical Engineering at the University of California, Los Angeles. Email:\texttt{\{adolfo,tabuada\}@ee.ucla.edu}, \texttt{http://www.ee.ucla.edu/$\sim$\{adolfo,tabuada\}} }
\end{flushleft}
% make the title area
\date{}}

\maketitle

\begin{abstract}
Feedback control laws have been traditionally implemented in a periodic fashion on digital hardware. Although periodicity simplifies the analysis of the mismatch between the control design and its digital implementation, it also leads to conservative usage of resources such as CPU utilization in the case of embedded control. We present a novel technique that abandons the periodicity assumption by using the current state of the plant to decide the next time instant in which the state should be measured, the control law computed, and the actuators updated. This technique, termed self-triggered control, is developed for two classes of nonlinear control systems, namely, state-dependent homogeneous systems and polynomial systems. The wide applicability of the proposed results is illustrated in two well known physical examples: a jet engine compressor and the rigid body.
%In this paper we address the problem of sampling for nonlinear systems. Traditionally feedback control laws are implemented in a periodic fashion, as it facilitates the analysis of the implementation effects. However, this approach leads to conservative implementations since the period has to be chosen based on the worst case scenario, so that the desired performance is guaranteed for all operation points. We present a novel technique that drops the periodicity assumption and instead decides the next sampling time based on the current state of the plant. Such technique is developed for two classes of nonlinear control systems, namely, state-dependent homogeneous and polynomial systems. The results are applied to two well known physical examples: a jet engine compressor and a rigid body.
%On the other hand, self-triggered techniques, since the inter-execution time is decided based on 
%Nowadays ppl go periodic. No good my friend, especially for nonlin. We address self-trig. for different class of nonlin. Examples illustrate.
\end{abstract}

\section{Introduction}

%$\mathcal{x}, \mathbb{x} \mathrm{x} \mathbf{x} \mathscr{x} \mathfrak{x} \mathit{x} \mathnormal{x}$
%$\mathcal{e}, \mathbb{e} \mathrm{e} \mathbf{e} \mathscr{e} \mathfrak{e} \mathit{e} \mathnormal{e}$
Feedback control laws are mainly implemented nowadays in digital platforms since microprocessors offer many advantages with respect to analog platforms. %\r{Although digital implementations can be described by discrete-time models, most of the available results in nonlinear control theory were developed for continuous-time systems}. 
Despite its advantages, digital implementations raise several design difficulties; one of the most relevant problems consists in determining how frequently the controller needs to be executed so that a desired performance is achieved. Traditionally, the controller is executed periodically every $T$ units of time, since this facilitates the analysis of the mismatch between the customary continuous design and the digital implementation. However, it seems unnatural to update the signals of interest in a periodic fashion, especially for nonlinear control systems: control tasks should be executed only when something significant happens in the process to be controlled. Moreover, the choice of the period $T$ is based on a worst-case scenario (to guarantee performance for all possible operation points), and hence the control task is executed at the same rate regardlessly of the state of the plant. Therefore, the periodicity assumption leads to inefficient implementations in terms of processor usage, communication bandwidth, energy,... For instance, in the context of embedded systems, the CPU time is shared between different tasks. Hence, executing the control task when the states barely change is a waste of computational resources. In the context of distributed control systems, not only the processor time is a scarce resource but also the available communication bandwidth. Thus communication should only occur when relevant information needs to be transmitted from the sensor to the controller and/or from the controller to the actuator. 

There is a vast literature about periodic implementations of control laws; still, there are many open points, especially for nonlinear systems. We refer to~\cite{hristu2005hna} for a good introduction to the subject for linear systems. Due to the lack of a deep understanding of the effects of digital implementations, ad hoc rules are commonly used to determine stabilizing periods (for instance, 20 times the time constant of the dominant pole~\cite{astrom90}). Most of the results relating stability and sampling periods for linear systems are based on the construction of an equivalent discrete-time model. Nonlinear systems, in general, cannot be discretized exactly in closed form. Hence, a common approach consists in finding approximate discrete-time models~\cite{nesic1999scs}, and then carry out the study and design for this set of equations~\cite{laila2005sdc}.

To overcome the drawbacks of the periodic paradigm, several researchers (\cite{arzen99}, \cite{astrom02}, \cite{sandee05}, \cite{tabuada07}) suggested the idea of event-triggered control. Under this paradigm, the controller execution is triggered according to the state of the plant. %Indeed, the control law should be executed only when something significant happens in the system. %One of the first studies of event-triggered control was carried in~\cite{arzen99}, where a simple PID controller is executed whenever the tracking error signal is large. 
The event-triggered technique reduces resource usage and provides a high degree of robustness (since the system is measured continuously). Unfortunately, in many cases it requires dedicated hardware to monitor the plant permanently, not available in most general purpose devices.

In this paper we propose to take advantage of the event-triggered technique without resorting to extra hardware. In most feedback laws, the state of the plant has to be measured (or estimated) to compute the next value of the controller; hence, this information could also be used to decide when the controller needs to be executed again. This technique is known as self-triggered control, since the controller decides its next execution time. This technique could be regarded as a way to introduce feedback in the triggering process, in comparison to the open loop structure in periodic implementations. A first attempt to explore self-trigger models for linear systems was developed in~\cite{velasco2003stt}, by discretizing the plant, and in~\cite{lemmon07} for linear ${\cal H}_{\infty}$ controllers. %In the context of nonlinear systems, to the best of our knowledge, the first results appeared in~\cite{tabuada08} where, under a homogeneity assumption, scaling laws for the inter-execution times were derived as a function of the state norm. 
No studies have been carried so far for nonlinear systems, where aperiodicity is expected to be even more beneficial since the behaviour of a nonlinear dynamical system highly depends on the operation point. 

The contributions of this paper are twofold:
\begin{enumerate}
\item Under self-trigger implementations, the execution times for a control law are defined by a simple self-trigger condition (since it has to be computed online) that depends on the dynamics of the system, the desired performance, and the current measurement of the state. We develop self-trigger conditions for two classes of systems, namely, state-dependent homogeneous systems and polynomial systems.
\item The developed self-trigger conditions are but a facet of the inherent tradeoff between the amount of resources alloted to a controller and the achieved performance. On our route to develop these conditions we develop the necessary theory required to further the understanding of the real-time requirements of control systems.
\end{enumerate}

The results in this paper can thus be seen as a contribution to the broad area of control under information constraints. %This topic represents one of the main issues of the problem of control design under resource constraints. 
Many authors have recently studied other aspects of this problem, such as space quantization~\cite{brockett2000qfs}, delays in network control systems~\cite{lian2003tdm} or minimum attention control design~\cite{brockett1997mac}, to cite a few. 
%The main contribution of the paper consists in deriving such formulas for two classes of nonlinear systems, namely, state-dependent homogeneous systems and polynomial systems. It represents as well an important step towards the understanding of the real-time requirements of control systems. To determine such formulas, we look at the mismatch between the current state and the sampled state. Similar results could be developed by looking at the mismatch of the input of the plant under analog and digital implementations. %Indeed, changes in the measured state might not imply significant changes in the controller output. 
%From a theoretical point of view both analysis are analogous; in this paper we will focus on the state-based approach rather than the input-based. \r{Regarding the time aspects of the real-time implementation, a control task is usually composed of three parts: measurement, computation of the control law and update of the actuator value. For ease of exposition, we will assume in the beginning that the measurement and the actuation occur instantaneously. Hence, every time we sample the current state the controller is computed without delay.} 

The rest of the paper is organized as follows. We first review in section III an event-trigger condition that guarantees stability under sample-and-hold implementations, previously studied in~\cite{tabuada07}. In section IV we define homogeneous control systems and we analyse the properties of their trajectories. Such properties can be exploited to derive a simple self-trigger condition based on the norm of the current state. This idea is generalized by introducing the notion of state-dependent homogeneity, %. These systems possess infinitesimal symmetries 
that allows us to establish a more general triggering condition based on the full state of the plant. %compute the next execution time for the control law. %Given the time for a point in the state space, we are able to infer the inter-execution times for scaled versions of such point. \textit{Indeed, there exist spatio-temporal relations for the inter-execution times that help us describe their evolution in the state space region}. 
%The aforementioned scaling law motivates us to derive self-triggered conditions that determine the next execution time in order to achieve desired performance. These results are exploited to derive similar scaling laws for polynomial systems.
Finally, in section V the previous results are utilized to develop self-trigger conditions for $\varphi$-related vector fields and, in particular, for polynomial control systems. To illustrate the wide applicability of the results we select two very different examples from the literature: the control of a jet engine compressor and the control of a rigid body. The appendix includes the proofs for all the results appearing in the paper.

\section{Notation and input-to-state stability}
\subsection{Notation}
\label{notation}
We shall use the notation $\vert x\vert$ to denote the Euclidean norm of an element $x\in\Re^n$. A continuous function \mbox{$\alpha:[0,a[\to \Re_0^+$}, $a>0$, is said to be of class $\mathcal{K}$ if it is strictly increasing and
$\alpha(0)=0$. It is said to be of class $\mathcal{K}_\infty$ if $a=\infty$ and $\alpha(r)\to \infty$ as $r\to\infty$. 

A function is said to be of class $C^{\infty}$ or smooth if it can be differentiated infinitely many times. All the objects in this paper are considered to be smooth unless otherwise stated. %Given a set $S \subseteq \R^n$ and a function \mbox{$f:\R^n \rightarrow \R^p$}, we write $f(S)$ to indicate that the function $f$ is applied to every point in the set $S$, that is $f(S)=\{f(s), s\in S\}$.

Given vector fields $X$ and $Y$ in an n-dimensional manifold $M$, we let [X,Y] denote their Lie product which, in local coordinates $x = (x_1,x_2,...,x_n)$, we take as $\frac{\partial{X}}{\partial{x}}(x)Y(x) - \frac{\partial{Y}}{\partial{x}}(x)X(x)$. Let $\varphi:M \rightarrow N$ be a map, and let $X$ and $Y$ be vector fields in manifolds $M$ and $N$, respectively. The differential of a map $\varphi$ is denoted by $T\varphi$, and whenever we need to emphasize the point $x$ where the differential is evaluated we write %the differential at a particular point by 
$T_x\varphi$. We call $X$ and $Y$ $\varphi$-related if the following holds:
$$ T \varphi \: X = Y \circ \varphi$$
If $\varphi$ has a smooth inverse, the pullback of the vector field Y under the mapping $\varphi$ can be defined as:
$$ \varphi^* Y \triangleq (T\varphi)^{-1} Y \circ \varphi$$ 
The notations $\psi(t,\cdot)$ and $\psi_t(\cdot)$ will be used interchangeably to denote a map $\psi: \R \times M \rightarrow M$. Finally, we use $\e^s$ to represent the exponential of $s\in\R$.

\subsection{Input-to-state stability}
\label{ISS_lyap}
We consider a control system:
\begin{equation}
\label{CS}
\dot{\mathrm{x}}=f(\mathrm{x},\mathrm{u}),\qquad \mathrm{x}(t)\in \Re^n,\,\, \mathrm{u}(t)\in \R^m, t\in\R_0^+
\end{equation}
%where $\mathrm{x}(t)$ denotes the state and $\mathrm{u}(t)$ the input. 
We use $x \in \R^n$ to denote the state of the control system, $\mathrm{x}(t)$ to denote a solution of~(\ref{CS}), and $\mathrm{u}(t)$ for the trajectory of the input. Whenever it is needed to emphasize the initial condition for the trajectory we write $\mathrm{x}(t,x_0)$, where $\mathrm{x}(0,x_0)=x_0$. We will not resort to the standard definition of input-to-state stable (ISS) (see~\cite{sontagISS} for details) in this paper but rather to the following characterization.
\begin{definition}
A smooth function $V:\Re^n\to\Re_0^+$ is said to be an ISS Lyapunov function for the system~(\ref{CS}) if there exist class $\mathcal{K}_\infty$ functions $\underline{\alpha}, \overline{\alpha}, \alpha$ and $\gamma$ satisfying:
\begin{eqnarray}
\label{CI}
&\underline{\alpha}(\vert x\vert)\le V(x)\le \overline{\alpha}(\vert x\vert)&\\
\label{CII}
&\frac{\partial V}{\partial x}f(x,u)\le -\alpha(\vert x\vert)+\gamma(\vert u\vert)&
\end{eqnarray}
The system~(\ref{CS}) is said to be ISS with respect to the input $\mathrm{u}$ if and only if there exists an ISS Lyapunov function for~(\ref{CS}). 
\end{definition}

\section{Event-triggered stabilization of nonlinear control systems}
\label{times_ev}

Even though most of the current controllers are implemented on digital platforms, the techniques used to design these controllers ignore the implementation effects, such as delays or space and time quantization. In this paper we focus on the time quantization aspects. We start analyzing the behaviour of the plant under the event-triggered implementation introduced in~\cite{tabuada07} and reviewed in this section. Consider the control system described in~(\ref{CS})
%Let our control system be described by:
%\begin{equation}
%\label{CS2}
%\dot{x}=f(x,u)
%\end{equation}
for which a feedback controller $u=k(x)$ has been designed. % to render the system ISS with respect to measurement errors.
The implementation of such feedback law on an embedded processor is typically done by sampling the state at time instants $t_i$, computing $\mathrm{u}(t_i)=k(\mathrm{x}(t_i))$ and updating the actuator values at time instants $t_i+\Delta_i$. For ease of exposition, it is assumed that $\Delta_i=0$ for all $i\in\N$, that is, every time the state is sampled the controller is computed without delay\footnote{The results can be generalized for nonzero $\Delta$ by following the procedure described in [Tab07].}. Furthermore, the sequence of times $\{t_i\}_{i\in\N}$ is typically periodic meaning that $t_{i+1}-t_i = T$ for any $i\in\N$, where $T>0$ is the period. In this paper the periodicity assumption is dropped in favour of self-trigger implementations. We would like to identify a class of aperiodic executions that guarantees stability and desired performance while reducing the number of executions. In other words, we would like to determine when it is indeed necessary to execute the control law in order to achieve the desired performance. To derive a stabilizing triggering rule, we look at the mismatch between the current state and the sampled state. Similar results could be developed by looking at the mismatch of the input of the plant under analog and digital implementations. %Indeed, changes in the measured state might not imply significant changes in the controller output. 
Since both analysis are analogous from a conceptual point of view, in this paper we just focus on the state-based approach for ease of exposition. %To derive a stabilizing execution rule, we study the inter-sample behaviour of the controlled system. 

We define the measurement error $\mathrm{e}$ as the difference between the last measured state and the current value of the state:
\begin{equation}
\label{Error}
\mathrm{e}(t)=\mathrm{x}(t_i)-\mathrm{x}(t) \qquad \text{for } t\in [t_i, t_{i+1}[
\end{equation}
With this definition, the closed loop $\dot{\mathrm{x}} = f(\mathrm{x},k(\mathrm{x}(t_i)))$ becomes:
\begin{equation}
\label{closed_loop}
\dot{\mathrm{x}}=f(\mathrm{x},k(\mathrm{x}+\mathrm{e}))
\end{equation}
Let the control law $u=k(x)$ render the system ISS with respect to the measurement error $e$. Under this assumption, there exists a Lyapunov function $V$ for the system that satisfies the following inequality:
\begin{equation}
\label{Vdot_nonlin}
\dot{V}\le -\alpha(\vert x\vert)+\gamma(\vert e\vert) 
\end{equation}
where $\alpha$ and $\gamma$ are ${\cal K}_\infty$ functions. Stability of the closed loop system~(\ref{closed_loop}) can be guaranteed if the error is restricted to satisfy:
\begin{equation}
\label{RuleNon}
\gamma (\vert e\vert) \le \sigma \alpha( \vert x\vert), \qquad \sigma > 0
\end{equation}
since the dynamics of $V$ will be bounded by:
\begin{equation}
\label{Vdot_nonlin2}
\dot{V}\le (\sigma-1) \alpha(\vert x\vert)
\end{equation}
guaranteeing that $V$ decreases provided that $\sigma<1$. If $\alpha^{-1}$ and $\gamma$ are Lipschitz continuous on compacts, inequality~(\ref{RuleNon}) is implied by the simpler inequality:
\begin{equation}
\label{Rule}
b \vert e\vert \le \sigma a \vert x\vert
\end{equation}
for $a$ and $b$ appropriately chosen according to the Lipschitz constants of $\alpha^{-1}$ and $\gamma$. Inequality~(\ref{Rule}) can be enforced by executing the control task whenever:
\begin{equation}
\label{ERule}
\vert e\vert = \sigma \frac{a}{b} \vert x\vert
\end{equation} 
Upon the execution of the control task, the state is measured and the error becomes 0, since \mbox{$\mathrm{x}(t_i)=\mathrm{x}(t)$} implies \mbox{$\mathrm{e}(t)=\mathrm{x}(t_i)-\mathrm{x}(t)=0$}. An event-triggered implementation based on this equality would require testing~(\ref{ERule}) frequently. Unless this testing process is implemented in hardware, one might run the risk of consuming the processor time testing~(\ref{ERule}). To overcome this drawback, we opt for self-trigger strategies, where the current state measurement $\mathrm{x}(t_i)$ is used to determine the next execution time $t_{i+1}$ for the control task. 

The inter-execution time implicitly defined by~(\ref{ERule}) is the time it takes for $\frac{\vert \mathrm{e}\vert}{\vert \mathrm{x}\vert}$ to evolve from\footnote{Recall that at the execution instant $t=t_i$ we have \mbox{$\mathrm{e}(t)=\mathrm{x}(t_i)-\mathrm{x}(t)=0$} and thus $\frac{\vert e \vert}{\vert x \vert}=0$.} $0$ to $\sigma\frac{a}{b}$. Since this time depends on the last sample of the system $x(t_i)$, it will be denoted as $\tau(x(t_i))$. In order to derive self-trigger conditions based on the sampled state, we study in the next sections the evolution of $\tau$ as a function of $x(t_i)$. This evolution can be described by analyzing the trajectories for different types of nonlinear systems. For simplicity of presentation, we start considering homogeneous systems of constant degree.

\section{Homogeneous control systems}
\label{class_homog}
Homogeneous vector fields are vector fields possessing a symmetry with respect to a family of dilations. They appear as local approximations for general nonlinear systems \cite{hermes1991nah} since an analytic function can always be decomposed in an infinite sum of homogeneous functions. Moreover, many physical systems can be described as homogeneous systems (see \cite{baillieul80} for some examples). 
\subsection{Definitions}
To define homogeneity we first review the notion of dilation.
\begin{definition}
\label{dilation}
Given an n-tuple \mbox{$r\!=\!(r_1,..,r_n) \in (\R^+_0)^n$}, the dilation map $\delta_{\lambda}^r: \R^n \rightarrow \R^n$ is defined by:
\begin{equation}
\delta_\lambda^r (x) = (\lambda^{r_1} x_1,...,\lambda^{r_n} x_n), \quad \lambda > 0
\end{equation}
A dilation is called standard if $r_i =r_j$ for all $i,j\in\{1,\ldots,n\}$. %A homogeneous ray is defined as the 1-parameter family of dilations $\{\delta^r_{\lambda}(x): \lambda > 0\}$. 

\end{definition}
\begin{definition}
A function $f: \R^n \rightarrow \R^n$ is called homogeneous of degree $\zeta$ if for all $\lambda > 0$, there exist \mbox{$r\!=\!(r_1,..,r_n) \in (\R^+_0)^n$} such that:
\begin{eqnarray}
\label{Dil_out}
& f (\delta_\lambda^r (x)) = \lambda^\zeta \delta_\lambda^r(f(x))
\end{eqnarray}
where $\zeta > -\min_i{r_i}$. 
\end{definition}
With this definition, we see that linear functions are homogeneous of degree $\zeta=0$ with respect to the standard dilation. Consider now a differential equation:
\begin{equation}
\label{DiffE}
\dot{\mathrm{x}} = f(\mathrm{x})
\end{equation} 
whose right-hand side is homogeneous of degree $\zeta$. Then, the solution $\phi(t,\phi_0)$ satisfies: 
\begin{equation}
\label{traj_homog_non_geom}
\phi(t,\delta^r_\lambda (x_0)) = \delta^r_{\lambda} \circ \phi (\lambda^\zeta t,x_0)
\end{equation}
For a proof see~\cite{kawski95} or the proof of Theorem~\ref{thm_homog_traj} in this paper that generalizes this fact. Equality~(\ref{traj_homog_non_geom}) lies at the heart of the results presented in this paper.

\subsection{Scaling laws for the inter-execution time of homogeneous systems}
\label{sect_scal_law_homog}
In this section we consider homogeneous control systems with respect to the standard dilation, since the general case can be reduced to this one, as explained in~\cite{grune2000hsf}. Using the commutative property of homogeneous systems expressed by equality~(\ref{traj_homog_non_geom}), we can derive the following scaling law for the inter-execution times under the event-triggered policy described in Section~\ref{times_ev}. 
\begin{theorem}
\label{ScaleHom}
Let $\dot{\mathrm{x}} = f(\mathrm{x},\mathrm{u})$ be a control system for which a feedback control law $u=k(x)$, rendering the closed loop homogeneous of degree $\zeta$ with respect to the standard dilation, has been designed. The inter-execution times $\tau:\R^n \rightarrow \R$ implicitly defined by the execution rule $\vert \mathrm{e} \vert =  c \vert \mathrm{x}\vert$ with $c > 0$ 
%and for any initial conditions lying on a homogeneous ray $\{\delta^r_\lambda(x_0): \lambda > 0\}$, 
scale according to:
\begin{equation}
\label{scal_eq}
\tau\big(\delta^r_{\lambda} (x)\big) = \lambda^{-\zeta} \tau(x), \qquad \lambda > 0
\end{equation}
where $x \in \R^n$ represents any point in the state space.
\end{theorem}

%\begin{proof}
%See appendix.
%\end{proof}
Similar results were obtained in~\cite{emre_tuna} but with a different objective (finite-time stability and finite escape time for homogeneous systems). This theorem relates the inter-execution time at a point $x$ with the inter-execution time at any other point lying along the same homogeneous ray. That is, once the time is known for just one initial condition $x$, we can infer the times for all initial conditions of the type $ \lambda x$, for any $\lambda > 0$. In the next section these ideas are extended by allowing the degree of homogeneity $\zeta$ to be a function of the state.

\subsection{State-dependent homogeneity and space-time dilations}
\label{state_dep_homog}
In order to generalize the classical notion of homogeneity presented in section~\ref{class_homog}, we resort to the coordinate-free geometric notion of homogeneity introduced in \cite{kawski95}. 
\begin{definition}
\label{geom_homog}
Let $D:M \rightarrow TM$ be a (dilation) vector field such that $\dot{\mathrm{x}} = -D(\mathrm{x})$ is globally asymptotically stable. A vector field $X: M \rightarrow TM$ is called homogeneous of degree $\zeta$ with respect to the vector field $D$ if it satisfies the following relation:
\begin{equation}
\label{Homog_vf}
[D,X] = \zeta X
\end{equation}
The trajectories $\psi: \R \times M \rightarrow M$ of $D$ are called the homogeneous rays of the system. 
\end{definition}

To recover the classical notion of homogeneity, we set: 
\begin{equation}
\label{std_dil_vf}
D = \sum_{i=1}^n r_i x_i \frac{\partial}{\partial x_i}, \qquad r_i > 0
\end{equation}
where $\big\{\frac{\partial}{\partial x_i}\big\}_{i=1,\ldots,n}$ is the canonical basis for the tangent bundle $TM$. By substituting this expression for $D$ in~(\ref{Homog_vf}), the classical notion of homogeneity (equation~(\ref{Dil_out})) can be recovered as explained in~\cite{arnol'd_ode}. Under this framework, we can see that the definition of section~\ref{class_homog} just corresponds to a particular choice of the vector field $D$. %This new definition allows us to consider different families of homogeneous rays. 
%We are mainly interested in the trajectories of the vector field $X$. 
In this geometric context,~(\ref{traj_homog_non_geom}) can be generalized to the following theorem.
\begin{theorem}
\label{thm_homog_traj}
Let $X$ and $D$ be vector fields on a manifold $M$, giving rise to flows $\phi:\R \times M \rightarrow M$ and $\psi:\R \times M \rightarrow M$  respectively. The vector field $X$ is homogeneous of degree $\zeta$ with respect to $D$ if and only if:
\begin{equation}
\label{Homog_traj}
\phi_t \circ \psi_s  = \psi_s \circ \phi_{\e^{\zeta s}t},  \qquad s,t \in \R
\end{equation}
\end{theorem}
%\begin{proof}
%See appendix.
%\end{proof}
This theorem implies that the flows of homogeneous vector fields commute in a particular way: applying the flow $\psi$ before the flow $\phi$ to a point $x\in M$ is equivalent to applying the flow $\phi$ before $\psi$ but with a scaling in time, given by $\e^{\zeta s}$, where $\zeta$ is the degree of homogeneity. For $D$ as defined in~(\ref{std_dil_vf}), equation~(\ref{Homog_traj}) simplifies to~(\ref{traj_homog_non_geom}), since the flow $\psi_s$ becomes the dilation map $\delta^r_\lambda$, for $\lambda=\e^s$. In this paper we generalize even further Definition~\ref{geom_homog} to the following state-dependent notion of homogeneity. 
\begin{definition}
Let $D:M \rightarrow TM$ be a (dilation) vector field such that $\dot{\mathrm{x}} = -D(\mathrm{x})$ is globally asymptotically stable. A vector field $X: M \rightarrow TM$ is called homogeneous with degree function \mbox{$\xi:M \rightarrow \R$} with respect to the vector field $D: M \rightarrow TM$ if it satisfies the following relation:
\begin{equation}
\label{Homog_gen}
[D,X] = \xi X
\end{equation}
\end{definition}
\begin{example}
\label{example_homog_def}
Let the vector fields $D$ and $X$ be:
\begin{eqnarray}
D & = & \alpha x_1 \frac{\partial}{\partial x_1} + \alpha x_2 \frac{\partial}{\partial x_2}, \qquad \alpha>0 \nonumber \\
X & = & (-x_1-x_1^3)\frac{\partial}{\partial x_1} + (-x_2-x_1^2 x_2) \frac{\partial}{\partial x_2}
\end{eqnarray}
The vector fields $D$ and $X$ do not satisfy (\ref{Homog_vf}) for any value of $\zeta$, but they satisfy~(\ref{Homog_gen}) for the following function $\xi$:
$$\xi(x) = \frac{3x_1^2+1}{x_1^2+1}$$
\end{example}
As for the standard homogeneity case, a relation between trajectories can also be derived. 
\begin{theorem}
\label{thm_gen_trajs_hom}
Let $X$ and $D$ be vector fields on a manifold $M$, giving rise to flows $\phi:\R \times M \rightarrow M$ and $\psi:\R \times M \rightarrow M$  respectively. The vector field $X$ is homogeneous with degree function $\xi: M \rightarrow \R$ with respect to $D$ if and only if:
\begin{equation}
\label{Homog_traj_gen}
\hspace{1.8cm} \phi_t \circ \psi_s  =  \psi_s \circ \phi_{\e^{\rho(s)}t}  \qquad s,t \in \R 
\end{equation}
with:
\begin{equation}
\hspace{-0.6cm} \qquad \rho(s) = \int_0^s \xi \circ \psi_{\tau} \: d\tau 
\end{equation}
%where $x\in M$.% and $\e^{\rho(s)}$ represents the dilation in time induced by the dilation in space $\psi$. 
\end{theorem}
%\begin{proof}
%See appendix.
%\end{proof}

The vector field $D$ generates a flow $\psi$ that can be considered as a spatial dilation operator acting on points in $M$. This operator determines how times are scaled in the flow of $\phi$. That is, given a point $x$ in a manifold $M$, applying this spatial operator $\psi$ to $x$ entails a scaling in time $\e^{\rho(s)}$ in the flow of a vector field $X$. 
\begin{remark}
In the case of the standard homogeneity, we have $\xi(x)=\zeta$ and thus equation~(\ref{Homog_traj_gen}) becomes~(\ref{Homog_traj}) since $\rho$ is linear in time, $\rho(s)=\zeta s$. Furthermore, for linear systems, the degree of homogeneity $\xi$ is $0$, therefore $\rho(s)=0$ and the vector fields $X$ and $D$ commute, that is:
\begin{eqnarray}
\label{linear_traj}
\phi_t \circ \psi_s = \psi_s \circ \phi_t
\end{eqnarray}
\end{remark}

\begin{example}
To illustrate the previous theorem, we recover Example~\ref{example_homog_def}. The flow for the dilation vector field $D$ is:
\begin{eqnarray}
%& &\frac{d}{dt} \phi(t,x) = D(\phi(t,x)), \hspace{1cm} \phi(0,x)=x\nonumber \\
& &\psi_s(x) = \e^{\alpha s} x
\end{eqnarray}
where $x=(x_1,x_2)^T$. Hence the flow of the vector field $X$ satisfies equation~(\ref{Homog_traj_gen}):
\begin{equation}
\hspace{1.8cm}\phi_t (\e^{\alpha s} x)  = \e^{\alpha s} \phi_{\e^{\rho(s)}t} (x) \\
\end{equation}
with:
\begin{equation}
\hspace{-0.8cm}\rho(s)  =  \frac{1}{\alpha} \log \bigg(\frac{\e^{3\alpha s}x^2_{1}+\e^{\alpha s}}{x^2_{1}+1}\bigg)
\end{equation}
We see that the dilation in space $\e^{\alpha s}$ induces a scaling in time $e^{\rho(s)}$ in the trajectory $\phi$.
\end{example}
Using these theorems, a self-trigger condition for state-dependent homogeneous systems can developed.
\subsection{Scaling laws for inter-execution times of state-dependent homogeneous systems}
\label{sect_scal_law_homog_gen}
For simplicity, we use the standard dilation vector field:
\begin{eqnarray}
\label{gen_dilation_vf}
D = \sum_{i=1}^n r x_i \frac{\partial}{\partial x_i}
\end{eqnarray}
The flow of this vector field (i.e., the homogeneous ray) is:
\begin{eqnarray}
%\frac{d}{ds} \psi_s(x) & = & D(\psi_s(x)), \hspace{1cm} \psi_0(x)=x\nonumber \\
\psi_i(s,x)  = \e^{r s} x_i \qquad \qquad i=1,\hdots,n \qquad x = (x_1,\hdots,x_n)^T
\end{eqnarray}
%since it satisfies:
%\begin{equation}
%\frac{d}{ds} \psi_s(x)  =  D(\psi_s(x)), \hspace{1cm} \psi_0(x)=x\\
%\end{equation}
%Moreover, we assume $r_1=r_2=\hdots=r_n$ (called standard homogeneity), since the general case of~(\ref{gen_dilation_vf}) can be reduced to this one by means of a change of coordinates, as explained in~\cite{grune2000hsf}. 
Using the commutative property~(\ref{Homog_traj_gen}) of homogeneous flows, we can derive a generalized scaling law for the inter-execution times.

\begin{theorem}
\label{ScaleHom_gen}
Let $\dot{\mathrm{x}} = f(\mathrm{x},\mathrm{u})$ be a control system for which a feedback control law $u=k(x)$ has been designed, rendering the closed loop homogeneous with degree function $\xi$ with respect to the standard dilation vector field. The inter-execution times $\tau:\R^n \rightarrow \R$ implicitly defined by the execution rule $\vert e \vert =  c \vert x\vert$, with $c > 0$, 
%and for any initial conditions lying on a homogeneous ray $\{\delta^r_\lambda(x_0): \lambda > 0\}$, 
scale according to:
\begin{equation}
\label{scal_eq_gen}
\hspace{1.8cm}\tau(\e^{rs} x) = \e^{-\rho(s)} \tau(x) \hspace{0.4cm} \forall s \in \R, r>0 \\
\end{equation}
with:
\begin{equation}
\label{rho_def}
\hspace{-0.8cm}\qquad  \rho(s) = \int_0^s \xi(\e^{r v}x) dv
\end{equation}
and where $x\in\R^n$ represents any point in the state space.
\end{theorem}
%\begin{proof}
%See appendix. 
%\end{proof}
%\begin{remark}
%\r{Many interesting conclusions can be drawn from this theorem. For instance, if a vector field $X$ in $\R^n$ is linear in some coordinates $\{x_i,\hdots,x_j\}$, the inter-execution times preserving stability do not depend on those coordinates. }
%\end{remark}
%This theorem relates the inter-execution time at a point with the inter-execution time at any other point lying along the same homogeneous ray. That is, once the time is known for just one initial condition $p$, we can infer the times for all initial conditions of the type $ \lambda p$, for any $\lambda \geq 0$. 
\begin{remark}
For linear systems, the degree function $\xi$ is 0, therefore $\rho(s)=0$ and $\e^{-\rho(s)}=1$. That is, for linear systems the inter-execution times remain constant as we move along homogeneous rays.
\end{remark}
\begin{remark}
When the degree function $\xi$ is a constant (Definition~\ref{geom_homog}), the scaling in time depends entirely on the value of $s$. %In this case, equation~(\ref{scal_eq}) simplifies to:
%\begin{eqnarray}
%\tau(\lambda^r p) &=& \lambda^{-m} \tau(p) \\
%\text{with }\qquad \lambda = e^s
%\end{eqnarray}
That is, to see how times scale it is enough to determine the position of the state $x$ on a homogeneous ray, but not on which particular ray $x$ lies, since time in all rays scale in the same way. In other words, the scaling in time is just a function of the norm of the state, and therefore~(\ref{scal_eq_gen}) simplifies to~(\ref{scal_eq}), for $\lambda=e^s$ and $\rho(s) = \zeta s$. On the other hand, for state-dependent homogeneity the scaling in time is determined by the ray where $x$ lies and the position of $x$ in that ray, that is, $\rho$ is a function of $s$ and $x$, since time in each homogeneous ray scales in a different manner.
\end{remark}

Theorem~\ref{ScaleHom_gen} allows us to use the estimate of the inter-execution times at some $x$ in order to determine the inter-execution times for the whole ray through $x$. Therefore, it is enough to find estimates of these times on any $n-1$ sphere, and then extend the estimates along homogeneous rays. Moreover, since a lower bound $\tau^*$ for the inter-execution times can be easily computed for linear systems\footnote{See~\cite{tabuada07} for one possible method to find such $\tau^*$.}, a $n-1$ sphere can always be chosen where a linear over-approximation for the control system can be obtained. To do so, we rewrite equation~(\ref{Homog_gen}) in local coordinates. A homogeneous function $g$ of degree $\xi(x)-1$ satisfies:
\begin{equation}
(\xi(x)+r_i) g(x) = \sum_{i=1}^n r_i x_i \frac{\partial g}{\partial x_i}
\end{equation}
%where $X^i$ represents the $i^{th}$ coordinate of the vector field $X$. 
Hence, for the closed loop system \mbox{$\dot{\mathrm{x}} = \tilde{f}(\mathrm{x},\mathrm{e}) = f(\mathrm{x},k(\mathrm{x}+\mathrm{e}))$} we can find a bound for $\vert \tilde{f}(x,e)\vert$ linear in $\vert x \vert$ and $\vert e \vert$:
\begin{align*}
\label{lin_approx}
\vert \tilde{f}(x,e) \vert  & =  \vert H(x,e) x + G(x,e) e \vert   \\
&\leq \vert H(x,e) \vert \vert x \vert + \vert G(x,e) \vert \vert e \vert  \\
&\leq \vert H(x^*_a,e^*_a) \vert \vert x \vert + \vert G(x^*_b,e^*_b) \vert \vert e \vert
\end{align*}
where:
\begin{eqnarray}
\label{lin_model}
& H(x,e):= \left[ \begin{array}{ccc} 
\frac{r_1}{\xi+r_1}\frac{\partial \tilde{f}_1}{\partial x_1} & \ldots &  \frac{r_1}{\xi+r_1}\frac{\partial \tilde{f}_1}{\partial x_n} \\
\vdots																						 & \ddots &  \vdots 																						\\
\frac{r_n}{\xi+r_n}\frac{\partial \tilde{f}_n}{\partial x_1} & \ldots &  \frac{r_n}{\xi+r_n}\frac{\partial \tilde{f}_n}{\partial x_n} \\
\end{array} \right]  \nonumber \\
& G(x,e):= \left[ \begin{array}{ccc} 
\frac{r_1}{\xi+r_1}\frac{\partial \tilde{f}_1}{\partial e_1} & \ldots &  \frac{r_1}{\xi+r_1}\frac{\partial \tilde{f}_1}{\partial e_n} \\
\vdots																						 & \ddots &  \vdots \\
\frac{r_n}{\xi+r_n}\frac{\partial \tilde{f}_n}{\partial e_1} & \ldots &  \frac{r_n}{\xi+r_n}\frac{\partial \tilde{f}_n}{\partial e_n} \\
\end{array} \right]  
\end{eqnarray}
and $(x^*_a,e^*_a)$ and $(x^*_b,e^*_b)$ are such that:
$$\vert H(x,e) \vert\leq \vert H(x^*_a,e^*_a) \vert, \qquad \vert G(x,e) \vert\leq \vert G(x^*_b,e^*_b) \vert$$
for all $(x,e)$ in a neighbourhood $\Omega$ around the origin. %our $n-1$ sphere. 
So given this set $\Omega$ we can find where the norm of these matrices $H$ and $G$ attain its maximum values and then work with the following linear model for the sole purpose of computing $\tau^*$:
\begin{equation}
\label{lin_eq}
\dot {\mathrm{x}} = H(x^*_a,e^*_a) \mathrm{x} + G(x^*_b,e^*_b) \mathrm{e}
\end{equation}
It is important to emphasize that we are not trying to find a linearized model, as it would not guarantee stability for the original nonlinear system. To summarize, the computation of the self-triggered execution strategy is made in 4 steps:

\begin{enumerate}

\item Define an invariant set $\Omega_x \subseteq \R^n$ around the equilibrium point, for instance a level set of the Lyapunov function. The execution rule~(\ref{ERule}) guarantees that there exists another invariant set $\Omega_e \subseteq \R^n$ for the error $e$, given by $\Omega_e = \{e\in \R^n : |e|\leq \sigma |x|, x \in \Omega_x\}$.

\item Compute $H$ and $G$ and calculate the points $(x^*_a, e^*_a)\in \Omega_x \times \Omega_e$ and $(x^*_b, e^*_b)\in \Omega_x \times \Omega_e$ where $\vert H \vert$ and $\vert G \vert$ are maximized. 

%\item Compute the inter-execution time $\tau^*$ for the linear model~(\ref{lin_eq}) as described in~(\ref{eq_tau}), by identifying $A+BK$ with $H(x^*_a,e^*_a)$ and $BK$ with $G(x^*_b,e^*_b)$. The time $\tau^*$ is a stabilizing sampling period of our original system for any initial condition lying in $\Omega$.
\item Find a stabilizing inter-execution time $\tau^*$ for the linear model~(\ref{lin_eq}). Among others, one possible technique to compute $\tau^*$ was described in~\cite{tabuada07}:

\begin{equation}
\label{eq_tau}
\tau^* = -\Psi -\frac{2}{i \Theta }\bigg(\arctan\Big(\frac{\alpha_1 + 2 \frac{a}{b} \sigma\alpha_2} {- i \Theta}\Big)\bigg)
\end{equation}

with:
\begin{eqnarray*}
\Psi     &=& -\frac{2}{i \Theta }\bigg(\arctan\Big(\frac{\alpha_1 + 2 \phi_0 \alpha_2} {- i \Theta}\Big)\bigg)\\
\Theta   &=& \sqrt{4\alpha_2 \alpha_0 -\alpha_1^2}\\
\alpha_0 &=& |H(x_a^*,e_a^*)|\\
\alpha_1 &=& |H(x_a^*,e_a^*)| + |G(x_b^*,e_b^*)|\\
\alpha_2 &=& |G(x_b^*,e_b^*)|
\end{eqnarray*}

This time $\tau^*$ is a stabilizing sampling period for the original system for any initial condition lying in $\Omega_x$ (since $\Omega_x$ is invariant).

\item Let $\Gamma$ be the largest ball inside $\Omega_x$, and let $d$ be its radius. Relate the current state $x$ with some point in the boundary\footnote{Note that the boundary of $\Gamma$ is an $n-1$ sphere.} of $\Gamma$ via a homogeneous ray, that is, find $s\in\R$ (the dilation in space) such that \mbox{$\e^{rs} (y) = x$} for some $y$ in the boundary of $\Gamma$. Since we are working with the standard dilation and since we have an estimate $\tau^*$ valid for any point in the boundary of $\Gamma$ the next execution time $\tau^{\downarrow}(x)$ of the control task can be computed by using~(\ref{scal_eq_gen}):
\begin{equation}
\label{final_eq}
\tau^{\downarrow}(x)=\e^{-\rho(s)}\tau^*	 
\end{equation}
\end{enumerate}
As $\tau^*$ can be precomputed offline, the evaluation of~(\ref{final_eq}) can be performed online in a very short time. The estimate $\tau^{\downarrow}(x)$ represents a lower bound for the inter-execution times implicitly defined by~(\ref{ERule}), and therefore it guarantees stability for the original system. It is important to notice that the self-trigger technique is trying to emulate the event-trigger condition defined in~(\ref{Rule}). In this context, the conservativeness of this approach (due to the mismatch between the self-trigger and the event-triggered policy) relies entirely on the accuracy of the period $\tau^*$guaranteeing stability for the linear system~(\ref{lin_eq}). That is, no conservativeness is added when the scaling law is applied. %This also implies that we could use any other available technique to find a lower bound $\tau^*$ (a stabilizing period in a set), and then nicely extend it to any point via homogeneous rays.

\subsection{Example: Jet engine compressor}
To illustrate the previous results, we consider the control of a jet engine compressor. The following model is borrowed from~\cite{krstic95}:
\begin{eqnarray}
\label{jet_model}
\dot{\mathrm{x}}_1 & = & - \mathrm{x_2} -\frac{3}{2}\mathrm{x_1}^2 -\frac{1}{2}\mathrm{x_1}^3 \nonumber \\
\dot{\mathrm{x}}_2 & = & \frac{1}{\beta^2}(\mathrm{x_1}-\mathrm{u})
\end{eqnarray}
where $\mathrm{x_1}$ is the mass flow, $\beta$ is a constant positive parameter, $\mathrm{x_2}$ is the pressure rise and $\mathrm{u}$ corresponds to the throttle mass flow, the input for this controller. In the model~(\ref{jet_model}) we have already translated the origin to the desired equilibrium point, hence the objective of the controller consists in steering both state variables $(x_1,x_2)$ to zero. The operation region is a ball of radius $5.4$ centered at the origin. A control law \mbox{$u = g (x_1,x_2)$} is designed to render the closed loop globally asymptotically stable:
\begin{equation}
\label{fdb_law_jet}
u = x_1 - \frac{\beta^2}{2}(x_1^2+1)(y+x_1^2 y + x_1 y^2) + 2 \beta^2 x_1
\end{equation}
where we have applied the nonlinear change of coordinates $y = 2\frac{x_1^2+x_2}{x_1^2+1}$. The closed loop equations become:
\begin{eqnarray}
\label{cl_jet}
\dot{\mathrm{x}}_1 & = & -\frac{1}{2} (\mathrm{x_1}^2+1) (\mathrm{x_1}+\mathrm{y}) \notag\\
\dot{\mathrm{y}}     & = & - (\mathrm{x_1}^2+1) \mathrm{y} 
\end{eqnarray}
An ISS Lyapunov function for this system can be found using SOStools~\cite{sostools}:
\begin{equation*}
V = 1.46 x_1^2-0.35 x_1 y+1.16 y^2
\end{equation*}
Bounds for the derivative of $V$ along trajectories (including measurement errors) are found using also SOStools:
\begin{equation*}
\dot V \leq -0.74 \cdot 10^8 |x|^4 + 0.90 \cdot 10^8 |x|^2 |e|^2
\end{equation*}
where $x=(x_1,y)^T$ and $e=(e_1,e_2)^T$. Hence stability can be guaranteed if we enforce the following inequality:
\begin{equation*}
0.90 |e|^2 \leq 0.74 \sigma^2|x|^2
\end{equation*}
The closed loop system~(\ref{cl_jet}) is homogeneous with respect to~(\ref{gen_dilation_vf}) for \mbox{$\xi(x_1,y) = \frac{2 x_1^2}{x_1^2+1}$}. For simplicity, we pick $r=1$ in the dilation vector field~(\ref{gen_dilation_vf}). %The degree function $m$ does not depend on the state variable $y$, since the vector field is linear in $y$. 
We select a value of $\sigma<1$ guaranteeing stability under rule~(\ref{Rule}), for instance $\sigma=0.33$. We find an estimate of the aperiodic time sequence following the 4 steps mentioned before.

\begin{itemize}

\item We define an invariant set $\Omega_x$ that encloses the operation region, \mbox{$\Omega_x=\{x\in\R^n|V(x)=27.04\}$}.

\item Calculate the weighted Jacobians as defined in~(\ref{lin_model}) and compute the points where the maximums are attained: $x_a^*=x_b^*=(12.5,0)$ and $e_a^*=e_b^*=(3.75,0)$.

\item Compute the inter-execution time for the linear model: \mbox{$\tau^*=7.63$ms} for the selected value of $\sigma$ and the desired operation region.

\item Finally, applying equation~(\ref{final_eq}), we obtain the following formula describing the inter-execution times for the control task:
\begin{equation}
\label{inter_exec_example}
%\tau\big(x(t_i)\big) = \frac{e^{3s}}{} \cdot \tau^* \\
\tau^{\downarrow}_{i+1}\big(\mathrm{x}_1(\tau^{\downarrow}_i),\mathrm{y}(\tau^{\downarrow}_i)\big) = \frac{29 \mathrm{x}_1(\tau^{\downarrow}_i)+d(\mathrm{x}(\tau^{\downarrow}_i))^2}{5.36d(\mathrm{x}(\tau^{\downarrow}_i)) \mathrm{x}_1(\tau^{\downarrow}_i)^2+d(\mathrm{x}(\tau^{\downarrow}_i))^2} \cdot \tau^*, \qquad i>0
\end{equation}
where $d(\mathrm{x}(\tau^{\downarrow}_i))$ is the norm of the previously measured state:\\ $$d(\mathrm{x}(\tau^{\downarrow}_i))=\sqrt{\mathrm{x}_1(\tau^{\downarrow}_i)^2+\mathrm{y}(\tau^{\downarrow}_i)^2}$$ 
\end{itemize}

\begin{figure}[ht]
   \begin{center}
 	\includegraphics[width=0.75\hsize]{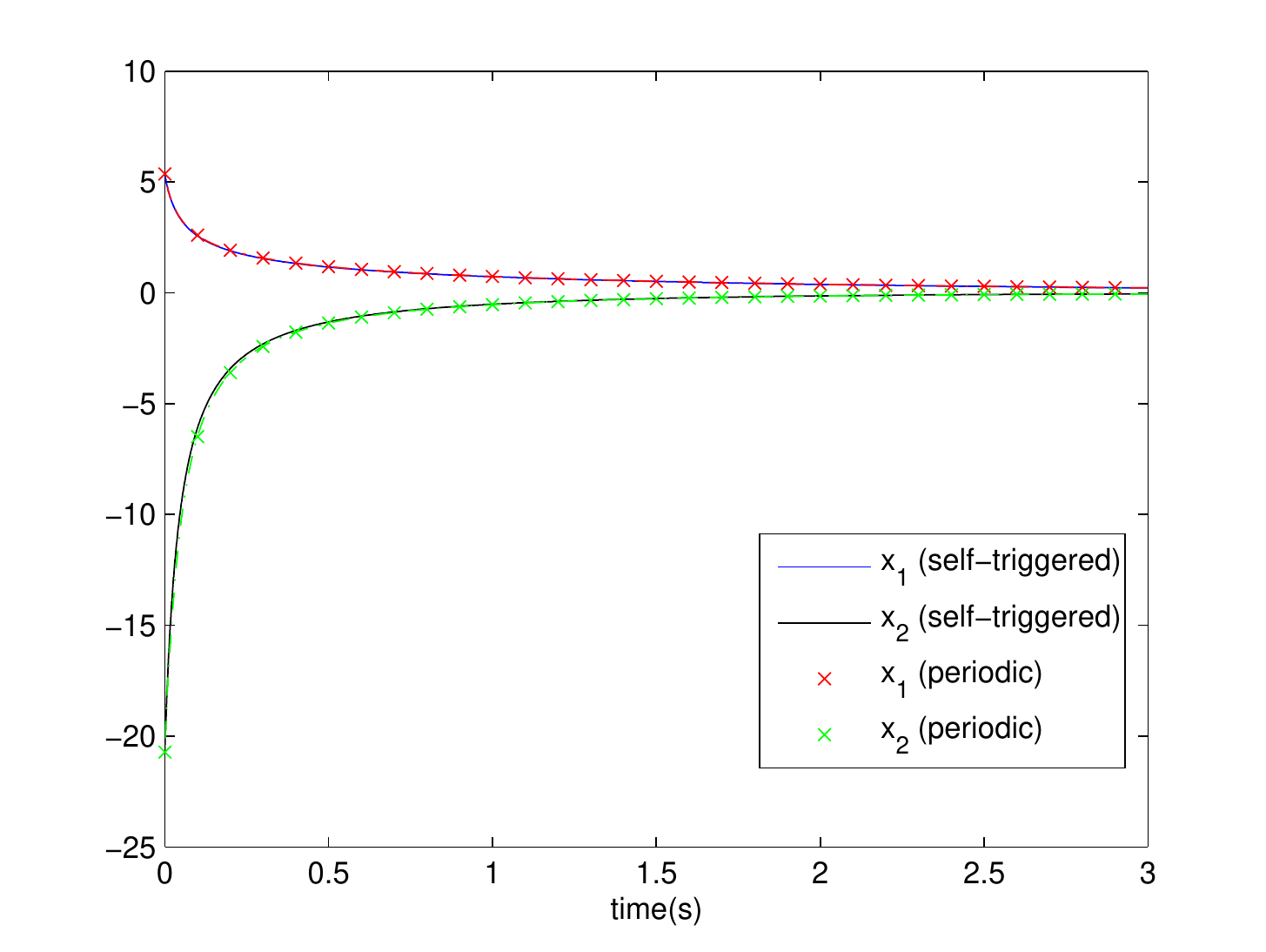}
   \end{center}
 \caption{Evolution of the states for self-trigger and periodic strategies.}
\label{evolution_orig_states}
\end{figure}

\begin{figure}[ht]
   \begin{center}
 	\includegraphics[width=0.75\hsize]{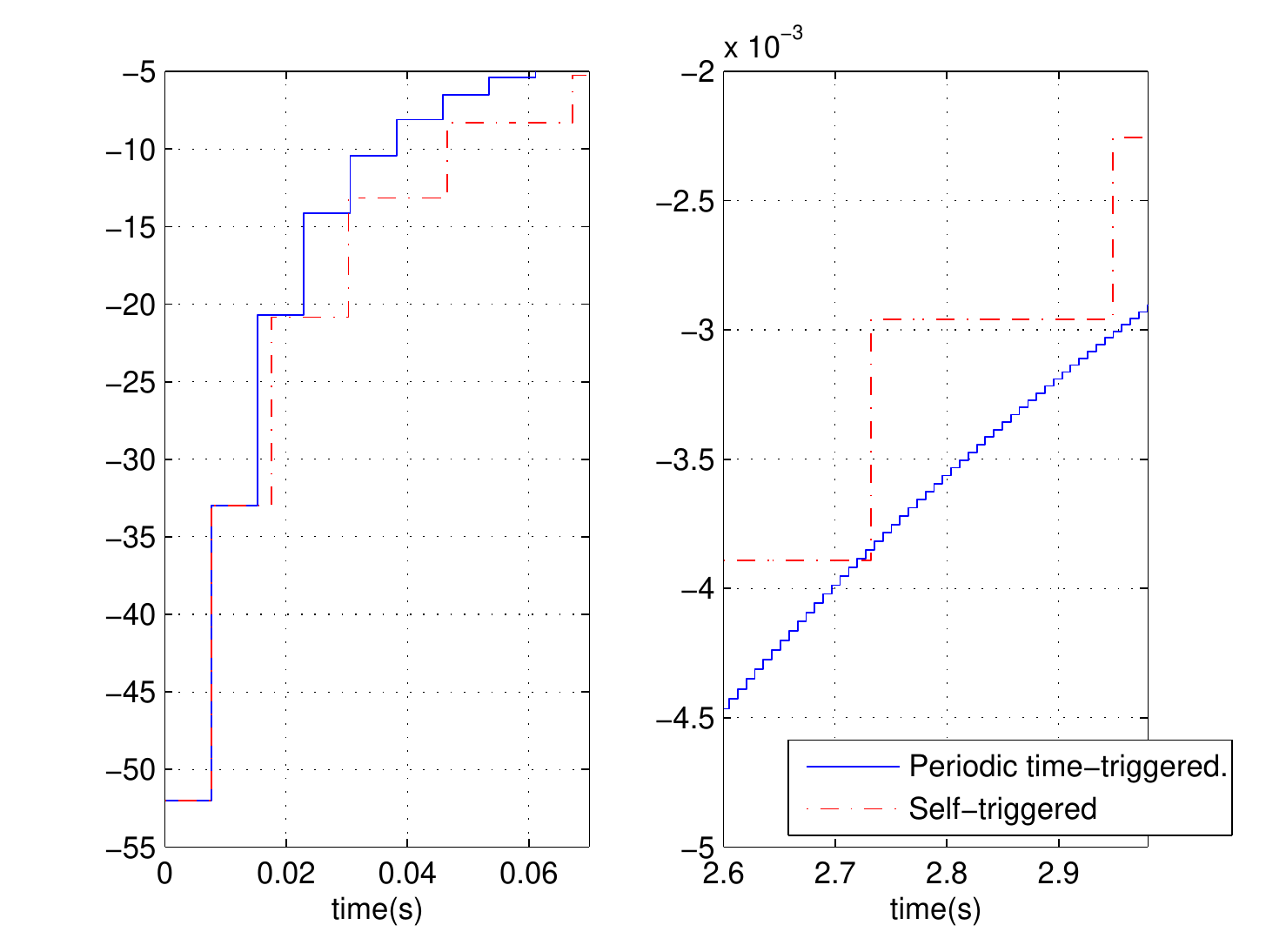}
   \end{center}
 \caption{Control input for self-trigger and periodic strategies.}
\label{fig_inputs}
\end{figure}

From this formula, we can see that times tend to enlarge as the system approaches the equilibrium point, since a linear term in $d(\mathrm{x}(\tau^{\downarrow}_i))$ only appears in the denominator of~(\ref{inter_exec_example}). In order to show the effectiveness of the approach, 50 different initial conditions were considered, equally distributed along the boundary of the operation region. In Figures~\ref{evolution_orig_states} and \ref{fig_inputs}, we compare the behaviour of both strategies, periodic and self-trigger. To choose a stabilizing period for the system we select the worst case inter-execution time obtained from~(\ref{inter_exec_example}). A different way to compute a stabilizing period for nonlinear systems appeared in~\cite{laila2003daa}; both techniques lead to similar values for the period. The systems exhibit a similar behaviour under both strategies for all initial conditions tested (see Figure~\ref{evolution_orig_states} for one particular initial condition, $x_0=(5.37,0.34)$). Figure~\ref{fig_inputs} shows the evolution of the input for the control system. At the beginning, both the periodic and self-trigger generate the same inter-execution times, but as the system tends to the equilibrium point the self-trigger policy increases the time  between executions, whereas the periodic policy keeps updating the controller at the same rate. The right side of Figure~\ref{fig_inputs} zooms the last part of the simulation, where the inter-execution times for the self-trigger strategy is already 24 times larger than the periodic. Hence the self-trigger implementation leads to a much smaller number of executions, while achieving a similar performance. The number of executions required for both implementations are shown in Table~\ref{comparison_table_jet}, for different values of $\sigma$ (and averaged over all initial conditions considered): the self-trigger policy reduces the number of executions by a factor of 8, for a simulation time of 3s. 
\begin{table}
\center
\begin{tabular}{cccc}
$\sigma$ & periodic & self-trigger\\
\hline\\
0.11 &  890 & 123\\
0.22 &  506 & 68\\
0.33 &  397 & 53\\
\hline\\
\end{tabular}
\caption{Number of executions of the control task for a simulation time of 3s.}
\label{comparison_table_jet}
\end{table}

%\begin{figure}[ht]
%   \begin{center}
% 	\includegraphics[width=0.8\hsize]{jet_evolution_orig_errors_no_noise.pdf}
%   \end{center}
% \caption{Comparison of trajectories under self-triggered and periodic strategies}
%\label{jet_evolution_orig_errors}
%\end{figure}

%Additionally, we would like to examine how accurately the procedure herein described approximates the inter-execution times defined by the event-triggered condition~(\ref{Rule}). Indeed, the least conservative inter-execution times are given by the event-triggered implementation, and our approach is just trying to emulate this technique in order to avoid extra hardware. Figure~\ref{scal_law1} compares the inter-execution times generated by the event-triggered and self-triggered techniques, for different initial conditions. Hence the mismatch between the self-triggered and the event-triggered policy depends on how tight is the bound $\tau^*$ for each point in the trajectory. Once again, we would like to emphasize that the difference between the two plots is only caused by the computation of the lower bound $\tau^*$ (equation~(\ref{eq_tau})), since no extra conservativeness is added by the scaling law~(\ref{inter_exec_example}).
\begin{figure}[ht]
   \begin{center}
 	\includegraphics[width=0.75\hsize]{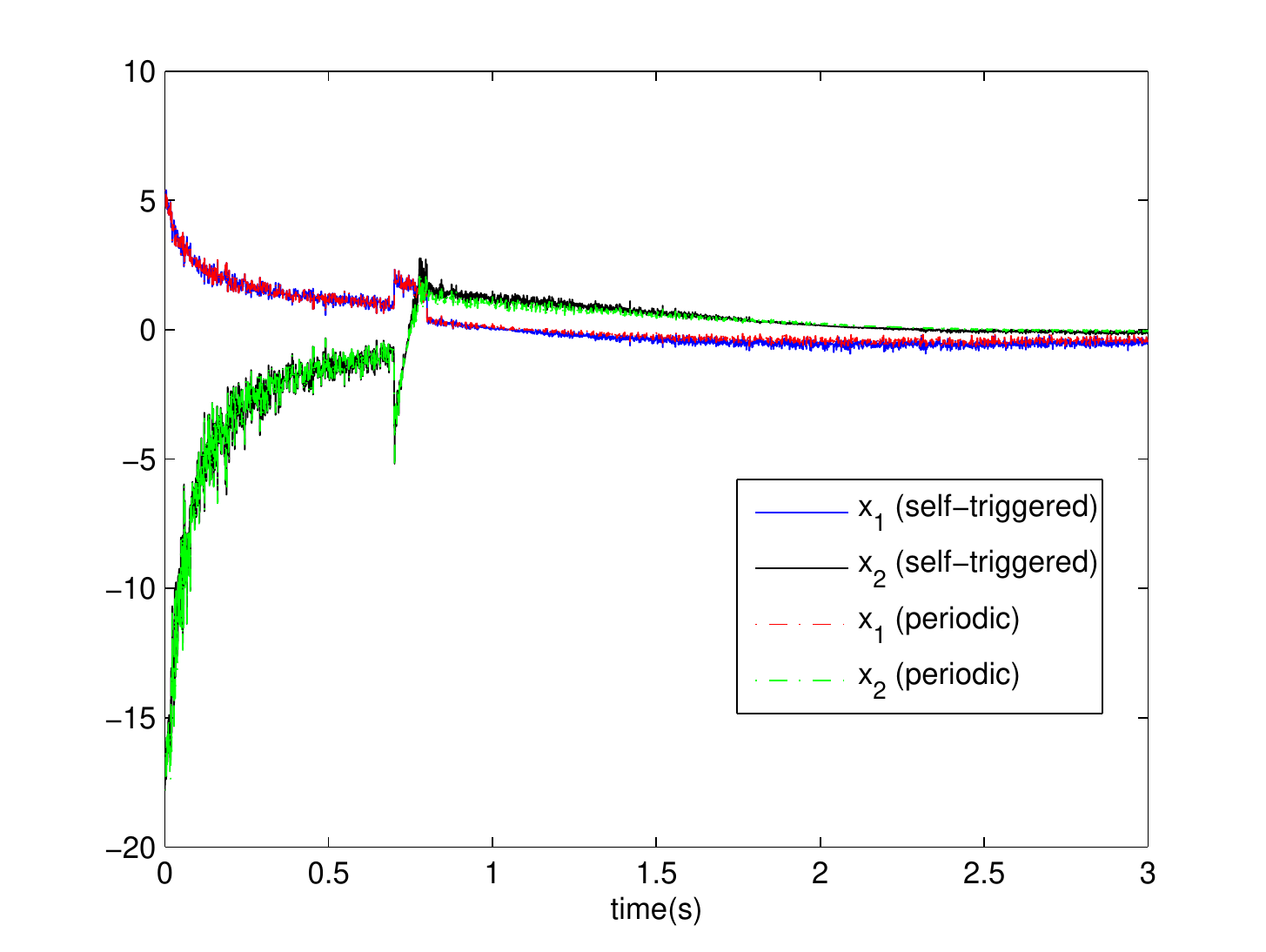}
   \end{center}
 \caption{Evolution of the states for self-trigger and periodic strategies in the presence of disturbances.}
\label{evolution_orig_states_noise}
\end{figure}

\begin{figure}[ht]
   \begin{center}
 	\includegraphics[width=0.75\hsize]{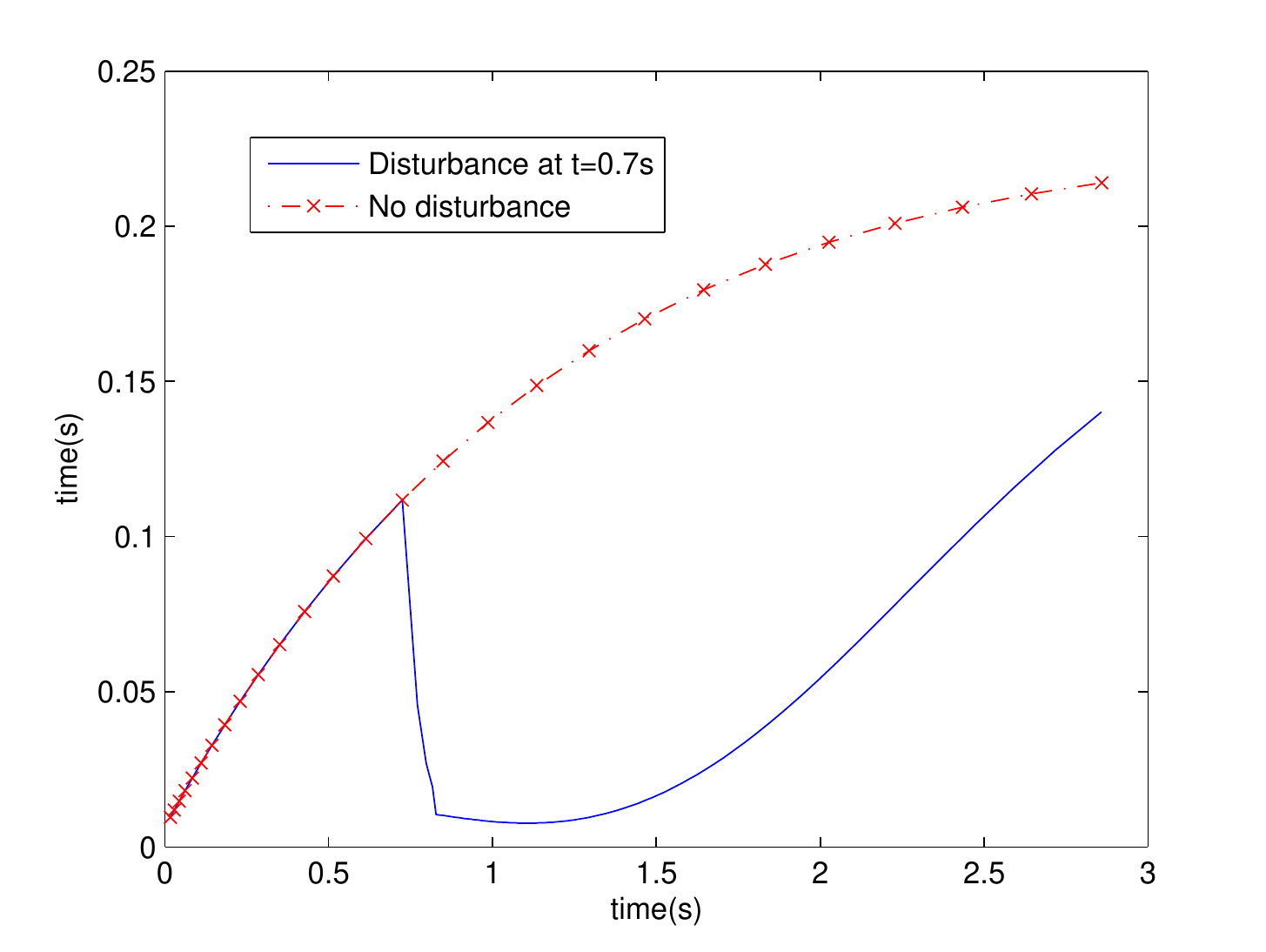}
   \end{center}
 \caption{Evolution of inter-execution times for the self-trigger strategy.}
 \label{scaling_law_pertb}
\end{figure}

We also verified through simulations that the self-trigger technique is robust with respect to disturbances, as expected because of the existence of an ISS Lyapunov function. Figure~\ref{evolution_orig_states_noise} shows the evolution of the states for both periodic and self-trigger strategies when sensor noise $n(t)$ is considered (noise power being 2\% of the signal power). In addition, a disturbance $d(t)$ at the actuator level is applied at $t=0.7s$. Again, the self-trigger strategy achieves a similar rate of decay with a much smaller number of executions. Figure~\ref{scaling_law_pertb} depicts the evolution of the inter-execution time under the self-trigger policy with and without the disturbance. When no disturbance $d(t)$ is applied, times enlarge as the states approach the origin, and the inter-execution times at $t=3$s rises above $210$ms. For the second case, the disturbance steers the system far from the origin at $t=0.7$s, and therefore the self-triggered task reduces the inter-execution times accordingly to guarantee the required performance. As the system reaches the equilibrium point, the inter-execution times start growing again since less executions are required to achieve the desired performance. Other possible real-time implementation effects that could degrade the performance of the system, such as jitter and input-output latency, are discussed in detail in~\cite{tabuada08_rtas}. The problem of real-time scheduling for a set of periodic, aperiodic and self-triggered tasks is also tackled in~\cite{tabuada08_rtas}. 

\section{Polynomial control systems}
\subsection{Scaling laws for inter-execution times of $\phi$-related systems}
In this section we show how the results developed for homogeneous systems can also be applied to other classes of vector fields.
In the case of homogeneous systems, the trajectories satisfy the scaling property~(\ref{Homog_traj_gen}). However, since the inter-execution times are determined by the evolution of the ratio $|\mathrm{e}|/|\mathrm{x}|$, it is sufficient to assume a scaling property just for the ratio, for instance:
\begin{equation}
\frac{|\mathrm{e}(t,\lambda x_0, \lambda e_0)|}{|\mathrm{x}(t,\lambda x_0,\lambda e_0)|}=\frac{|\mathrm{e}(\theta(x_0,e_0)t,x_0,e_0)|}{|\mathrm{x}(\theta(x_0,e_0)t,x_0,e_0)|}
\end{equation}
for some function $\theta:\R^n \times \R^n \rightarrow \R$ and $\lambda \in \R$. This observation can be exploited as follows. Let our open loop control system be:
\begin{equation}
\label{poly_ol}
\dot{\mathrm{x}}=f(\mathrm{x},\mathrm{u})
\end{equation}
A controller of the form $u = k (x)$ is designed to render the system ISS with respect to measurement errors. Due to the digital implementation, the feedback control law becomes $u = k (x+e)$ and the closed loop is:
\begin{eqnarray}
\label{poly_cl}
\dot{\mathrm{x}} & = & f(\mathrm{x},k(\mathrm{x}+\mathrm{e})) = \hat f(\mathrm{x}, \mathrm{e})
\end{eqnarray}
Since $\dot{\mathrm{e}} = -\dot{\mathrm{x}}$, the measurement error can be included in the state space representation of the system:
\begin{eqnarray}
\label{cl_gen}
\left[ \begin{array}{ccc} \dot{\mathrm{x}} \\ \dot{\mathrm{e}} \end{array}\right]= Z(\mathrm{x},\mathrm{e}) = \left[ \begin{array}{ccc} \phantom{-}\hat f(\mathrm{x},\mathrm{e}) \\ -\hat f(\mathrm{x},\mathrm{e}) \end{array}\right]
\end{eqnarray}

Let the solution to~(\ref{cl_gen}) be denoted by $\mathrm{z}$, and let $z = (x,e)$ be the state of this system. Now we define $\pi_x$ as the projection of $\mathrm{z}$ of the first $n$ coordinates ($\pi_x \circ \mathrm{z} = \mathrm{x}$) and $\pi_e$ as the projection of the last $n$ coordinates ($\pi_e \circ \mathrm{z}=\mathrm{e})$. We also define an output map $\eta_Z := |\pi_e|/|\pi_x|$. Then, the execution rule $|\mathrm{e}|=c|\mathrm{x}|$ becomes $\eta_Z \circ \mathrm{z}=c$. We thus look for another vector field $Y$ and an output map $\eta_Y$ that describe the evolution of times for the original system, that is $\eta_Z(\mathrm{z(t)})=\eta_Y(\mathrm{y(t)})$. These ideas can be summarized in the following theorem.

%\begin{displaymath}
%\xymatrix{
%\R & M \ar[l]_{\frac{|\pi_e|}{|\pi_x|}} \ar[r]^{Z} \ar[d]_{\phi} & TM \ar[d]^{T\phi}\\
%& N \ar[ul]^{\gamma_{xe}} \ar[r]^{Y} & TN}
%\R^n \ar[r]^{|\cdot|} \ar[d]_{\phi} & \R \\
%\R^{n+1} \ar[ur]_{|\cdot|}}
%\end{displaymath}

\begin{theorem}
\label{thm_phi_related}
Consider the vector field $Z:M\rightarrow TM$, as defined in~(\ref{cl_gen}), and a vector field \mbox{$Y:N \rightarrow TN$} giving rise to flows $\mathrm{z}$ and $\mathrm{y}$, respectively. If there exist maps $\varphi:M \rightarrow N$ and $\eta_Y:N \rightarrow \R$ such that the following diagram commutes:

\begin{equation}
\label{comm_diag}
\begin{diagram}
\R  & \lTo^{\eta_Z}         & M & \rTo^Z         & TM\\
    & \luTo_{\eta_Y}        & \dTo^{\varphi}     &   & \dTo_{T\varphi}\\  
    &  & N & \rTo^Y & TN
\end{diagram}
\end{equation}
that is:
\begin{enumerate}
	\item $T\varphi Z = Y \circ \varphi$;
	\item $\eta_Y \circ \varphi = \eta_Z$,
\end{enumerate}
then the inter-execution times $\tau_z(z)$ implicitly defined by $\eta_Z \circ \mathrm{z}(t,z)= c$ (with $c>0$) and the inter-execution times $\tau_y(\varphi(z))$ implicitly defined by $ \eta_Y \circ \mathrm{y} (t,\varphi(z)) = c$ coincide:
\begin{equation}
\tau_z (z) = \tau_y (\varphi(z))
\end{equation}
\end{theorem}
The auxiliary vector field $Y$ is $\varphi$-related to $Z$, and it faithfully represents the evolution of the ratio $|\pi_e|/|\pi_x|$ of the original system. Moreover, if $Y$ is homogeneous, scaling properties can be inferred for the ratio $|\pi_e|/|\pi_x|$, and thus we can derive scaling laws for the inter-execution times of the original vector field $Z$. In the next section we describe an application of these abstract ideas to the particular case of polynomial control systems. %the flow $\xi$ possesses symmetries that allow us to derive self-triggered conditions: %If $\xi$ is the flow of $Y$, then it satisfies~(\ref{traj_homog_non_geom}) and there exists a scaling law for the original system:
%\begin{eqnarray}
%\frac{|\pi_e|}{|\pi_x|}\circ \zeta (t,\lambda x, \lambda e) = \gamma_{xe} \circ \xi(t,\phi(\lambda x, \lambda e)) = \gamma_{xe} \circ \xi (\theta(x)t,\phi(x,e)) = \frac{|\pi_e|}{|\pi_x|} \circ \zeta(\theta(x)t,x,e)
%\end{eqnarray}
%Hence we could derive a scaling law for the inter-execution times for~(\ref{cl_gen}):
%\begin{eqnarray}
%\tau_z (\lambda x_0) = \frac{1}{\theta(x_0)}\tau_z(x_0)
%\end{eqnarray}

\subsection{Scaling laws for inter-execution times of polynomial systems}
\label{sec_homogenization}
As a concrete application of Theorem~(\ref{thm_phi_related}), we consider the class of polynomial systems. As explained for example in~\cite{baillieul80}, any polynomial vector field can be rendered homogeneous by introducing another state variable $w \in \R$, satisfying $\dot{\mathrm{w}}=0$. For instance, the polynomial system:
\begin{eqnarray}
\label{orig_sys_example}
\dot{\mathrm{x}}_1& = & \mathrm{x}_1 \mathrm{x}_2 + \mathrm{x}_2\nonumber \\
\dot{\mathrm{x}}_2& = & \mathrm{x}_1
\end{eqnarray}
can be easily rendered homogeneous of degree $1$ with respect to the standard dilation by adding the auxiliary variable $w$:
\begin{eqnarray}
\label{aux_sys_example}
\dot{\mathrm{x}}_1& = & \mathrm{x}_1\mathrm{x}_2+ \mathrm{x}_2 \mathrm{w}\nonumber \\
\dot{\mathrm{x}}_2& = & \mathrm{x}_1 \mathrm{w} \nonumber \\
\dot{\mathrm{w}}& = & 0 
\end{eqnarray}
Trajectories for the original system~(\ref{orig_sys_example}) can be recovered from the trajectories of the auxiliary system~(\ref{aux_sys_example}) provided that $\mathrm{w}(0)=1$. This homogenization procedure can be formalized as follows. The closed loop~(\ref{poly_cl}) is assumed to be polynomial, that is:
\begin{eqnarray}
\label{poly}
& \hat f_i(x, e)=\sum_{j=1}^{p}\alpha_{ij} x_1^{r_{1_{ij}}} \hdots x_n^{r_{n_{ij}}} e_1^{s_{1_{ij}}} \hdots e_n^{s_{n_{ij}}} & \\
& p \in \N_{\geq 1} \qquad r_{k_{ij}} \in \N_0, \qquad s_{k_{ij}} \in \N_0, \qquad \forall k=1,..,n  \qquad i=1,\ldots,n & \notag
\end{eqnarray}
We define $l$ as the highest degree in any of the monomials of $f_i$ for any $i=1,..,n$. Each monomial of degree $m<l$ is multiplied by $w^{m-l}$. With this new state variable $w$ we need to define a dummy error variable $e_w$, since under the framework of section~\ref{class_homog} it is necessary to have as many errors as states. %This new state variable $w$ comes together with a dummy measurement error variable $e_w$. 
Thus each component $\hat f_i$ becomes a homogeneous polynomial of degree $l-1$ in the variables $x_1,..,x_n,w,e_1,..,e_n,e_w$. Hence the state space representation of the extended system is:
\begin{eqnarray}
\label{CS_homogenized}
\dot{\mathrm{x}} &=& \tilde f(\mathrm{x}, \mathrm{w}, \mathrm{e},\mathrm{e}_w)\notag \\
\dot{\mathrm{w}} &=& 0
\end{eqnarray}
with 
\begin{eqnarray}
\label{poly_hom}
& \tilde{f_i}(x, w, e, e_w)=\sum_{j=1}^{p}\alpha_{ij} x_1^{r_{1_{ij}}} \hdots x_n^{r_{n_{ij}}} w ^{l- m_{ij}} e_1^{s_{1_{ij}}} \hdots e_n^{s_{n_{ij}}} \qquad i=1,\ldots,n &  \notag \\
\label{defs_degrees}
&  m_{ij}=\sum_{k=1}^n(r_{k_{ij}}+s_{k_{ij}}) \qquad l =\max_{i,j} m_{ij}& 
\end{eqnarray}
and $\mathrm{w}(0)=1$ so that $\tilde f(x,w,e,e_w) = \hat f(x,e)$. Using this procedure, the system~(\ref{CS_homogenized}) is homogeneous with respect to the standard dilation. In the context of the previous section, we can see that the vector field:
\begin{eqnarray}
\label{f_rel_vf}
Y(x,w,e,e_w) = \left[ \begin{array}{ccc} \phantom{-}\tilde f(x,w,e,e_w) \\ 0\\-\tilde f(x,w,e,e_w) \\0\end{array}\right]
\end{eqnarray}
is $\varphi$-related to the original polynomial vector field, with $\varphi(x,e) = (x,1,e,0)$ since:
%\begin{equation}
%T\phi \: Z (x,e) = Y \circ \phi (x,e)\\
%\end{equation}
%since:
\begin{eqnarray}
T\varphi Z (x,e) & = & \left[ \begin{array}{c} \phantom{-}\hat f(x,e) \\ 0 \\ -\hat f(x,e) \\0 \end{array}\right] 
= \left[ \begin{array}{c} \phantom{-}\tilde f(x,1,e,0) \\ 0  \\ -\tilde f(x,1,e,0) \\ 0\end{array}\right] \notag\\
&=&  Y(x,1,e,0) = Y \circ \varphi (x,e)
\end{eqnarray}
To apply Theorem~\ref{thm_phi_related} we define the output map $\eta_Y (x,w,e,e_w) = \frac{|e|}{|x|}$ for the extended system $Y$. %For this particular case, the auxiliary vector field $Y$ describes not only the ratio $|e|/|x|$ but the complete flow of the original vector field $Z$. 
Under this embedding, the inter-execution times for the original system~(\ref{poly_cl}) can be computed using Theorem~\ref{thm_phi_related}, as defined in the following corollary.
\begin{corollary}
\label{cor_poly}
Let $\tau(x)$ be the inter-execution times for the system~(\ref{cl_gen}) implicitly defined by $|\mathrm{e}|=c|\mathrm{x}|$ for a point $x\in\R^n$, and let $\tilde{\tau}(x,w)$ be the inter-execution times for the system~(\ref{CS_homogenized}) under the same execution rule, for a point $(x,w)\in\R^{n+1}$. Then, $\tau$ and $\tilde \tau$ are related according to:%satisfies:
\begin{equation}
\label{eq_cor_poly}
%\tau(x) = \lambda^{-l+1} \:\:\tilde{\tau} (x/\lambda,1/\lambda), \qquad \forall \lambda > 0
\tau(x) = \lambda^{l-1} \:\:\tilde{\tau} (\lambda x,\lambda), \qquad \forall \lambda > 0
\end{equation}
with $l$ as defined in~(\ref{defs_degrees}).
\end{corollary}
%\begin{proof}
%See appendix. 
%\end{proof}
The underlying idea consists in embedding the original system in a higher dimensional space, where symmetries for the flows  can be established and that were not present in the original $n$-dimensional space. This corollary does not provide any scaling law for the original system~(\ref{poly_cl}), as it was the case in sections~\ref{sect_scal_law_homog} and~\ref{sect_scal_law_homog_gen}, since we always have to look at the auxiliary system~(\ref{CS_homogenized}) to compute $\tau$. In order to infer a self-trigger condition from this theorem, the same steps as in Section~\ref{sect_scal_law_homog_gen} can be followed, but now carrying the computations for a $n+1$ dimensional space. Even though the extended auxiliary system no longer needs to be stable, existence of invariant sets (required in Section~\ref{sect_scal_law_homog_gen} to compute $\tau^*$) can be established.%proved, as shown in the following proposition.

\begin{proposition}
\label{prop_inv}
Consider the dynamical system~(\ref{cl_gen}), with $\hat f$ as defined in~(\ref{poly}). If $\Sigma$ is an invariant set for~(\ref{cl_gen}), then for any $\lambda>0$ the set $\tilde \Sigma_{\lambda} = \{(\lambda x,\lambda,\lambda e,0)\in \R^{2n+2}: (x,e)\in\Sigma\}$ is an invariant set for the extended system~(\ref{f_rel_vf}).	
\end{proposition}

%Since the union of invariant sets is also invariant, 
This proposition provides us a way to construct an invariant set where the computation of $\tau^*$ can be performed. The estimate $\tau^*$ of the inter-execution times can then be extended via homogeneous rays using Corollary~\ref{cor_poly}. It is important to notice that, unlike for the case of homogeneous systems, it is not possible to guarantee that the self-trigger condition always results in longer inter-execution times than the periodic approach, since $\tau^*$ is computed in a higher dimensional space. Nonetheless, the technique herein explained outperforms the periodic approach for all tested examples, one of them being discussed in the following section.

\subsection{Example: Rigid body}
\begin{figure}[ht]
   \begin{center}
 	\includegraphics[width=0.75\hsize]{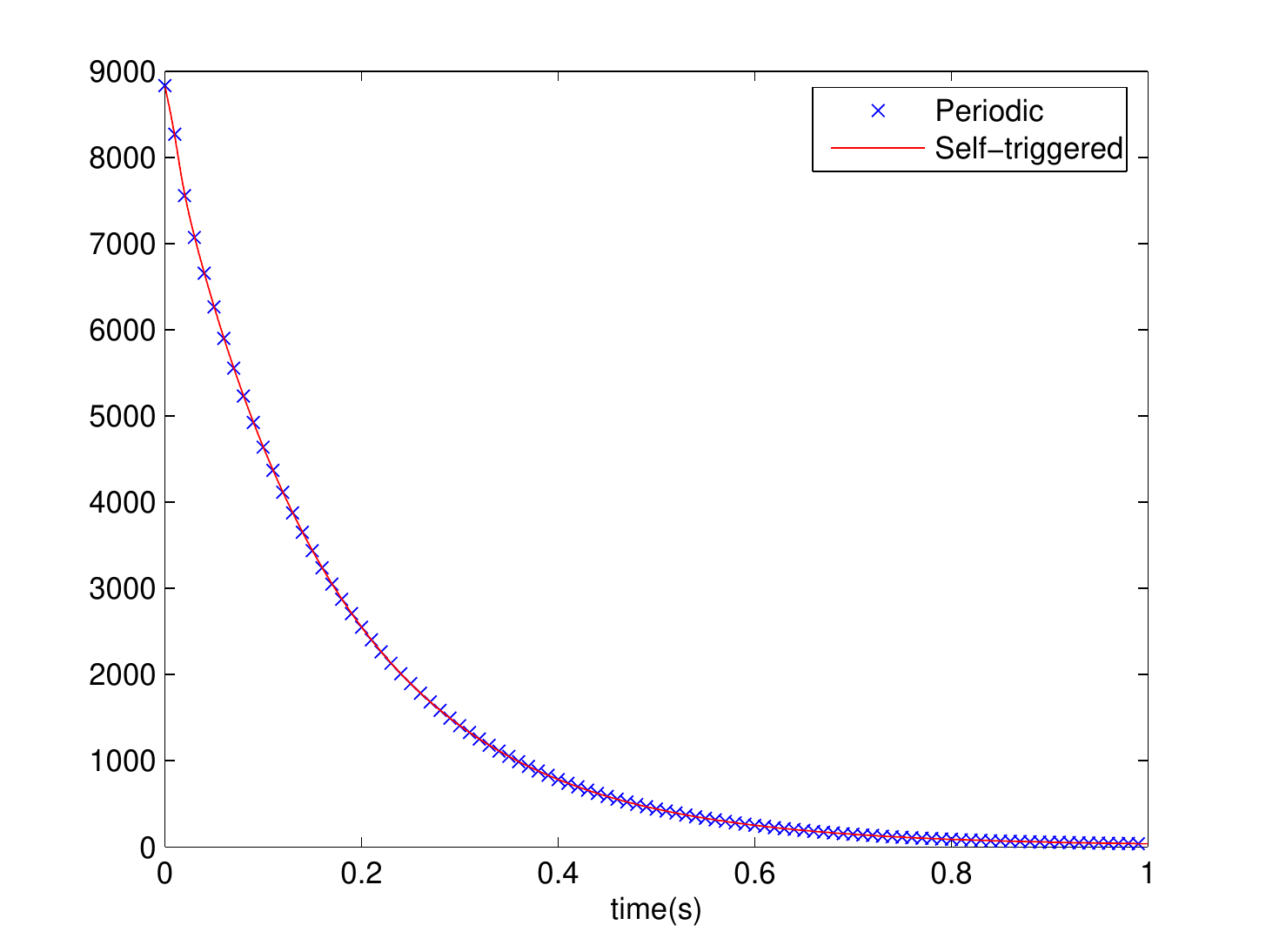}
   \end{center}
 \caption{Evolution of the Lyapunov function for self-trigger and periodic strategies.}
\label{satellite_lyap}
\end{figure}

\begin{figure}[ht]
   \begin{center}
 	\includegraphics[width=0.75\hsize]{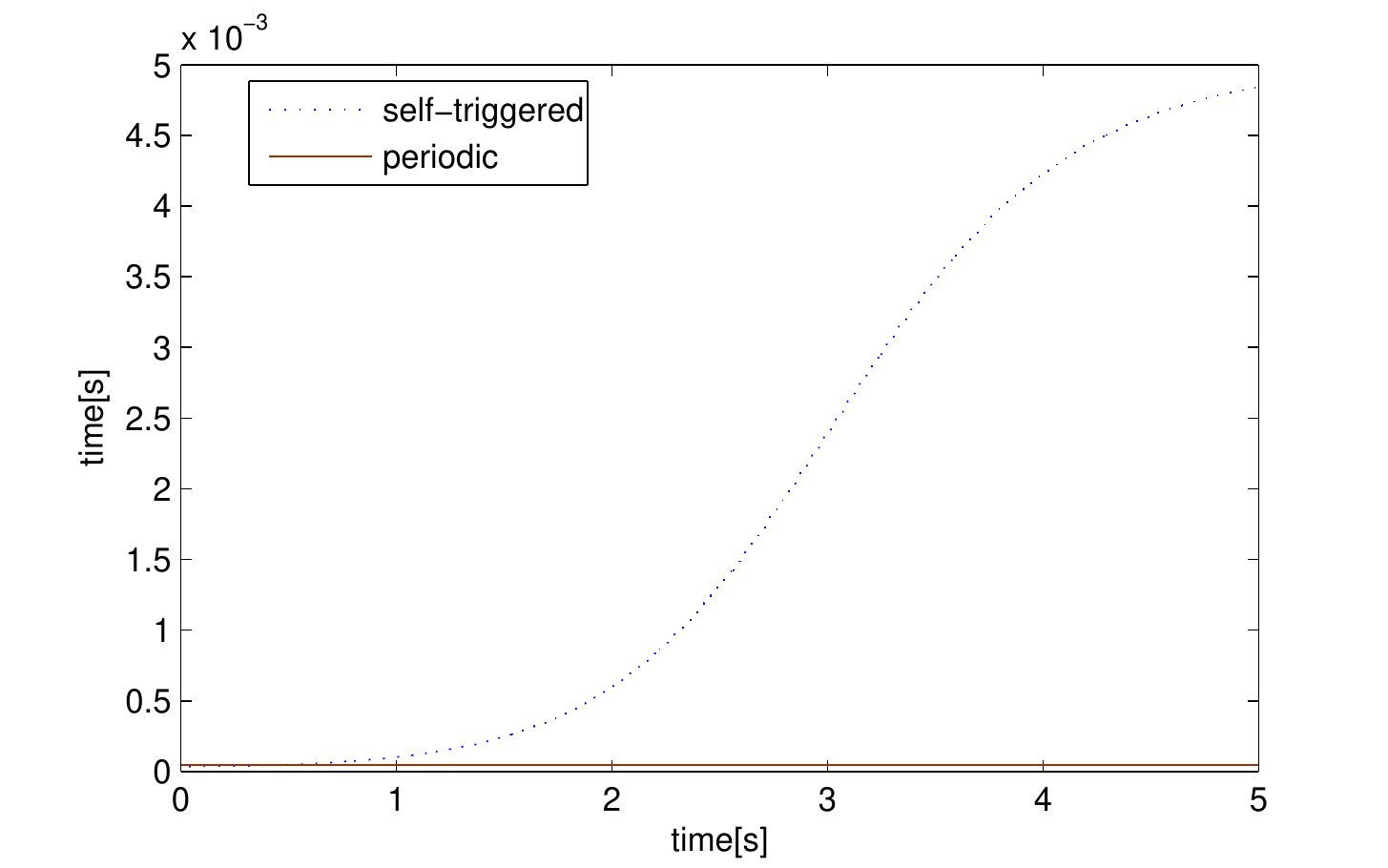}
   \end{center}
 \caption{Evolution of the inter-execution times for self-trigger and periodic strategies.}
\label{satellite_times}
\end{figure}

We apply the previous techniques to the control of the angular velocity for a rigid body. This example is borrowed from~\cite{byrnes1989nra}. After a preliminary feedback and normalization with respect to the moments of inertia, the state space representation of such system with two inputs can be simplified to the form:
\begin{eqnarray}
\dot{\mathrm{x}}_1& = & \mathrm{u}_1 \nonumber \\
\dot{\mathrm{x}}_2& = & \mathrm{u}_2 \nonumber \\
\dot{\mathrm{x}}_3& = & \mathrm{x}_1 \mathrm{x}_2
\end{eqnarray}
A nonlinear feedback law is designed in~\cite{byrnes1989nra} to render the system globally asymptotically stable:
\begin{eqnarray}
u_1 &=& -x_1 x_2 - 2 x_2 x_3  -x_1 - x_3 \nonumber \\
u_2 &=& 2 x_1 x_2 x_3 + 3 x_3^2 -x_2
\end{eqnarray}
as proved by the Lyapunov function:
\begin{eqnarray}
V(x_1,x_2,x_3) = \frac{1}{2}(x_1+x_3)^2+\frac{1}{2}(x_2-x_3^2)^2+x_3^2
\end{eqnarray}
The derivative of $V$ along the trajectories of the closed loop system (including measurement errors) is bounded by:
\begin{eqnarray}
\dot V \leq -91446 |x|^4 + 147190 |e|^2|x|^2
\end{eqnarray}
where $x=(x_1,x_2,x_3)^T$ and $e=(e_1,e_2,e_3)^T$. Hence stability will be guaranteed if:
\begin{eqnarray}
147190 |e|^2 \leq  91446 \sigma^2 |x|^2
\end{eqnarray}
for $\sigma<1$. For this particular example we select $\sigma=0.01$. The operation region is a ball around the origin of radius 15. We define an invariant set that encloses the operation region:
$$\Omega=\{(x,e)\in\R^{2n}:V(x) \leq 25650, |e| \leq \sigma |x|\}$$
and compute a stabilizing period for the system using the technique described  in~\cite{tabuada07} (equation~(\ref{eq_tau})). For the selected value of $\sigma$ and for this operation region we obtain $\tau^*=4.5 \cdot 10^{-5}$s. To homogenize the original polynomial system, the auxiliary variable $w$ is included in the state space representation for the closed loop:
\begin{eqnarray}
\label{example_homogenized}
\dot{\mathrm{x}}_1&=& -(\mathrm{x}_1+\mathrm{e}_1) (\mathrm{x}_2+\mathrm{e}_2) \mathrm{w} - 2 (\mathrm{x}_2 + \mathrm{e}_2) (\mathrm{x}_3+\mathrm{e}_3)\mathrm{w}  -(\mathrm{x}_1+\mathrm{e}_1)\mathrm{w}^2 \notag\\
&\phantom{=}& - (\mathrm{x}_3+\mathrm{e}_3) \mathrm{w}^2\nonumber \\
\dot{\mathrm{x}}_2&=& 2 (\mathrm{x}_1+\mathrm{e}_1) (\mathrm{x}_2+\mathrm{e}_2) (\mathrm{x}_3+\mathrm{e}_3) + 3 (\mathrm{x}_3+\mathrm{e}_3)^2\mathrm{w} - (\mathrm{x}_2 + \mathrm{e}_2)\mathrm{w}^2     \nonumber \\
\dot{\mathrm{x}}_3&=& \mathrm{x}_1 \mathrm{x}_2 \mathrm{w} \nonumber \\
\dot{\mathrm{w}}  &=& 0
\end{eqnarray}
%To compute $\tau^*$ for this example we define an invariant set that encloses the operation region, namely \mbox{$\Omega=\{(x,e)\in\R^{2n}:V(x) = 25650$, |e| \leq \sigma |x|\}}. 
This extended system is now homogeneous of degree $2$ with respect to the standard dilation. In order to derive a self-trigger condition based on~(\ref{eq_cor_poly}), it is necessary to find a lower bound $\tilde \tau^*$ for the inter-execution times $\tilde \tau$ of the extended system. An invariant set $\tilde \Sigma$ covering the whole operation region can be defined for the extended system using Proposition~\ref{prop_inv}:
$$\tilde \Sigma=\{(\lambda x,\lambda,\lambda e,0)\in\R^{2n+2}: (x,e)\in\Sigma_{\lambda}, 0 \leq \lambda \leq 0.066\}$$
where $\Sigma_{\lambda}=\{(x,e)\in\R^{2n}:V(x) \leq \frac{3}{2}(1-\lambda^2)+\frac{1}{2}(1-\lambda^2)^2, |e| \leq \sigma |x|\}$.
%where the computation of $\tilde \tau^*$ is performed. 
As before, the weighted Jacobians for the extended system~(\ref{example_homogenized}) are computed, and the maximums are attained at: $$x_a^*=(1.3,2.3,-1.4),\qquad w_a^*=0.066, \qquad e_a^*=(0.031,0.007,-0.035).$$ 
$$x_b^*=(0.7,2.2,-1.3),\qquad w_b^*=0.066, \qquad e_b^*=(0.028,0.005,-0.047).$$
%To compute $$ An invariant set can be found according to proposition.
Finally, applying the results in Corollary~\ref{cor_poly}, the following self-trigger condition for the control law is found:
\begin{eqnarray}
\label{jet_engine_scaling_law}
\tau^{\downarrow}_{i+1}(x(\tau^{\downarrow}_i)) = \lambda^{-2} \tilde \tau^*, \qquad \lambda = \sqrt{1+|x(\tau^{\downarrow}_i)|^2} 
\end{eqnarray}
where $\tilde \tau^*=5.1$ms. %for the selected value of $\sigma$ and the desired operation region. %The invariant set is not a ball, but to simplify online computations we pick the largest ball inside the invariant. 
Equation~(\ref{jet_engine_scaling_law}) shows that times tend to enlarge as the system approaches the origin. The rigid body system was simulated under the self-trigger and periodic implementations. Sensor noise was also included in the simulations (noise power being 2\% of the signal power). We compare the behaviour of the Lyapunov functions for the periodic and self-trigger implementations in Figure~\ref{satellite_lyap}. Again, both techniques achieve a similar performance. The evolution of times is displayed in Figure~\ref{satellite_times}. At the beginning, the periodic strategy generates a larger inter-execution time than the self-trigger, but as the system approaches the origin, the self-trigger policy increases the inter-execution times according to~(\ref{jet_engine_scaling_law}). For a simulation time of $5$s, the self-trigger policy reduces the total number of executions by a factor of $8$. %To avoid this initial overhead in the inter-execution times, the microprocessor should always select the , since both policies guarantee stability.
%As before, table of number of executions + evolution states.. or Lyapunov function.

%In this case, since the system is 3-d is better to plot the Lyapunov function rather than the evolution of the states

%\section{Future work?}
%How to use results for RT, NCS, energy saving... 
%Actuation errors?
%Conclusions?
%In this paper we have established . The results were applied to two physical systems, the jet engine compressor and the classical rigid body with two inputs. Simulations show how the number of executions is dramatically reduced. These ideas could be applied in the context of real-time systems, in order to save, or in the context of networked control systems, where...

\section{Appendix}
In this appendix we present the proofs of the results stated in this paper. Theorem~\ref{ScaleHom} and Theorem~\ref{thm_homog_traj} are not proved since they are special cases of Theorem~\ref{ScaleHom_gen} and Theorem~\ref{thm_gen_trajs_hom}, respectively, for the case when the degree of homogeneity is constant.

%\vspace{0.4cm}

\hspace{0.3cm} \textit{Theorem 4.8}: Let $X$ and $D$ be vector fields on a manifold $M$, giving rise to flows \mbox{$\phi:\R \times M \rightarrow M$ and $\psi:\R \times M \rightarrow M$}  respectively. The vector field $X$ is homogeneous with degree function $\xi: M \rightarrow \R$ with respect to $D$ if and only if:
\begin{equation}
\label{Homog_traj_app}
\hspace{1.8cm} \phi_t \circ \psi_s  =  \psi_s \circ \phi_{\e^{\rho(s)}t}  \qquad s,t \in \R 
\end{equation}
with:
\begin{equation}
\label{rho_def_app_ini}
\hspace{-0.6cm} \qquad \rho(s) = \int_0^s \xi \circ \psi_{\tau} \: d\tau 
\end{equation}

\begin{proof}
Let $x$ be a point in the manifold $M$. Both the infinitesimal and the integral part need to be proved. We start with the infinitesimal part. 
\begin{eqnarray}
%\phi_t \circ \psi_s  (x) & = & \psi_s \circ \phi_{\e^{\rho(s)}t} \Rightarrow [D,X](x) = m(x)X(x)\nonumber \\
\label{Homog_traj_proof}
\phi_t \circ \psi_s (x)& = & \psi_s \circ \phi_{\mathbf{e}^{\rho(s)}t} (x)= \psi_s \circ \phi_{(\e^{\rho(s)}-1)t} \circ \phi_t (x)\nonumber \\
\Leftrightarrow \phi_t \circ \psi_s \circ \phi_{-t} (x)& = & \psi_s \circ \phi_{(\e^{\rho(s)}-1)t} (x)
\end{eqnarray}
We differentiate first the left hand side of the previous equation with respect to $s$:
\begin{eqnarray}
\frac{d}{ds} \phi_t \circ \psi_s \circ \phi_{-t} (x) & = & T_{\psi_s \circ \phi_{-t} (x)} \phi_t \: \frac{d}{ds}\psi_s (\phi_{-t}(x)) \notag\\
&=& T_{\psi_s \circ \phi_{-t} (x)} \phi_t \: D \circ \psi_s (\phi_{-t}(x))
\end{eqnarray}
where the last equality holds since $\psi$ satisfies $\dot \psi = D (\psi)$. Evaluating the previous expression at $s=0$ we obtain:
\begin{eqnarray}
\label{thm1_lhs_nec}
\frac{d}{ds}\bigg|_{s=0} \phi_t \circ \psi_s \circ \phi_{-t}(x) = T_{\phi_{-t} (x)} \phi_t D (\phi_{-t}(x))
\end{eqnarray}
Differentiating the right hand side of equation~(\ref{Homog_traj_proof}) with respect to $s$ we obtain:
\begin{eqnarray}
\label{rhs1}
\frac{d}{ds} \psi_s \circ \phi_{(\e^{\rho(s)}-1)t}(x) & = & \frac{\partial}{\partial s} \psi_s (\phi_{(\e^{\rho(s)}-1)t}(x))\notag\\
&&+ T_{\phi_{(\e^{\rho(s)}-1)t}(x)} \psi_s \:  \frac{d}{ds}\phi_{(\e^{\rho(s)}-1)t}(x) \nonumber \\
& = & D \circ \psi_s \circ \phi_{(\e^{\rho(s)}-1)t}(x) \notag\\
&&+ T_{\phi_{(\e^{\rho(s)}-1)t}(x)} \psi_s \: X \circ \phi_{(\e^{\rho(s)}-1)t}(x)\:\e^{\rho(s)}t\frac{d}{ds}\rho(s)\nonumber\\
\end{eqnarray}
Using equation~(\ref{rho_def_app_ini}), we find the derivative of $\rho(s)$ with respect to $s$:
\begin{eqnarray}
\frac{d}{ds}\rho(s) & = & \frac{d}{ds} \bigg(\int_0^s \xi \circ \psi_{\tau} (x) d\tau \bigg)= \xi \circ \psi_s (x)
\end{eqnarray}
Since equation~(\ref{Homog_traj_app}) has to be true for any value of $s$, it also needs to hold for $s=0$, hence $\rho(0)=0$. Evaluating expression~(\ref{rhs1}) at $s=0$ we obtain:
\begin{eqnarray}
\label{thm1_rhs_nec}
& &\frac{d}{ds}\bigg|_{s=0} \psi_s \circ \phi_{(\e^{\rho(s)}-1)t}(x) = D(x) + X(x) t \xi(x)
\end{eqnarray}
We now equate the right hand sides of~(\ref{thm1_lhs_nec}) and~(\ref{thm1_rhs_nec}), differentiate them with respect to $t$ and finally evaluate them at $t=0$:
\begin{eqnarray}
T_{\phi_{-t} (x)} \phi_t D \circ \phi_{-t}(x) & = & D(x) +X(x) t \xi(x) \nonumber \\
\Rightarrow \frac{d}{dt}\bigg|_{t=0} T_{\phi_{-t} (x)} \phi_t D \circ \phi_{-t}(x) & = & \frac{d}{dt}\bigg|_{t=0} D(x) +X(x) t \xi(x) \nonumber \\
\Rightarrow \![D, X](x) & =& \xi(x)X(x)
\end{eqnarray}
For the integral part, we need to check the other direction of the statement:
\begin{eqnarray}
\label{int_impl}
[D, X](x) = \xi(x) X (x)\Rightarrow \phi_t \circ \psi_s (x)= \psi_s \circ \phi_{\e^{\rho(s)}t} (x)
\end{eqnarray}
We will integrate the first equality of~(\ref{int_impl}) with respect to $t$ and $s$. First we review the definition of pullback of a vector field $X$ under the mapping $\psi$ (see chapter 2):
\begin{equation}
\psi_s^* X(x) =  T_x \psi_{-s} X \circ \psi_s(x)\\
\end{equation}
%Then, we notice the following???-> Lie derivative: d/ds of pushfwd for s=0 is equal to the Lie derivative along that vector field
We define:
\begin{eqnarray}
\gamma(s,x) \triangleq \psi_s^* X(x)
\end{eqnarray}
and note that:
\begin{eqnarray}
\gamma(0,x) = \psi_0^* X(x) = X(x)
\end{eqnarray}
Moreover, the function $\gamma$ satisfies the differential equation:
\begin{eqnarray}
\frac{d}{ds} \gamma(s,x) &=& \frac{d}{ds} \psi_s^* X(x) = \psi_s^* [D,X](x) = \psi_s^* (\xi X) (x) \nonumber \\
&=& \psi_s^* \xi (x)\psi_s^*X (x)= \xi \circ \psi_s (x) \gamma(s,x)
\end{eqnarray}
whose solution is:
\begin{eqnarray}
\gamma(s,x) = \e^{\int_0^s \xi\circ \psi_{\tau}(x) d\tau} \gamma(0,x) = \e^{\rho(s)} \gamma(0,x) = \e^{\rho(s)} X(x)
\end{eqnarray}
according to the definition in~(\ref{rho_def_app}). Hence:
\begin{eqnarray}
\label{pullback_X}
\psi_s^*X (x) = \e^{\rho(s)}X (x)
\end{eqnarray}
The next step is to prove the following implication:
$$ \psi_s^*X(x) = \e^{\rho(s)}X(x) \Rightarrow \phi_t \circ \psi_s (x) = \psi_s \circ \phi_{\e^{\rho(s)}t} (x)$$
By definition of pullback, we can rewrite equation~(\ref{pullback_X}) as:
\begin{eqnarray}
(T_x \psi_s)^{-1} X \circ \psi_s(x)  & =  & \e^{\rho(s)}X (x)  \\
\Leftrightarrow X \circ \psi_s(x) & = & T_x \psi_s \: \e^{\rho(s)} X (x) 
\end{eqnarray}
Again, we define a function:
\begin{eqnarray}
\chi(t,x) \triangleq \psi_s \circ \phi_{\e^{\rho(s)}t} (x)
\end{eqnarray}
and show that it satisfies the following differential equation:
\begin{eqnarray}
\frac{d}{dt}\chi (t,x) & = & T_{\phi_{\e^{\rho(s)}t}(x)} \psi_s X \circ \phi_{\e^{\rho(s)}t}(x) \e^{\rho(s)} =  T_{\phi_{\e^{\rho(s)}t}(x)} \psi_s \e^{\rho(s)} X \circ \phi_{\e^{\rho(s)}t}(x) \nonumber\\
&=& X \circ \psi_s \circ \phi_{\e^{\rho(s)}t}(x)= X (\chi(t,x))
\end{eqnarray}
So $\chi$ is the flow of the vector field $X$. Since the initial conditions of $\chi$ and $\phi$ have to be the same ($\chi(0,x) = \psi_s (x)$) and the vector field $X$ generates a unique flow, equation~(\ref{Homog_traj_app}) is finally obtained:
\begin{equation*}
\phi_t \circ \psi_s (x) = \psi_s \circ \phi_{\e^{\rho(s)}t} (x) 
\end{equation*} %\qed
\end{proof}

\hspace{0.3cm} \textit{Theorem 4.11}: Let $\dot{\mathrm{x}} = f(\mathrm{x},\mathrm{u})$ be a control system for which a feedback control law $u=k(x)$ has been designed, rendering the closed loop homogeneous with degree function $\xi$ with respect to the standard dilation vector field. The inter-execution times $\tau:\R^n \rightarrow \R$ implicitly defined by the execution rule $\vert e \vert =  c \vert x\vert$, with $c > 0$, 
%and for any initial conditions lying on a homogeneous ray $\{\delta^r_\lambda(x_0): \lambda > 0\}$, 
scale according to:
\begin{equation}
\label{scal_eq_appendix}
\hspace{1.8cm}\tau(\e^{rs} x) = \e^{-\rho(s)} \tau(x) \hspace{0.4cm} \forall s \in \R, r>0 \\
\end{equation}
with:
\begin{equation}
\label{rho_def_app}
\hspace{-0.8cm}\qquad  \rho(s) = \int_0^s \xi(\e^{r v}x) dv
\end{equation}
and where $x\in\R^n$ represents any point in the state space.

\begin{proof}
Under the control law $u=k(x)$, the closed loop becomes:
$$\dot{\mathrm{x}} = f(\mathrm{x},k(\mathrm{x}))$$
Now we include the measurement errors as variables in the state space representation. Since $\dot e = -\dot x$ (see~(\ref{Error})), the closed loop  vector field can be written as:
\begin{eqnarray}
\label{cl_homog}
\left[ \begin{array}{ccc} \dot{\mathrm{x}} \\ \dot{\mathrm{e}} \end{array}\right]= Z(\mathrm{x},\mathrm{e}) = \left[ \begin{array}{ccc} \phantom{-}f(\mathrm{x},k(\mathrm{x}+\mathrm{e})) \\ -f(\mathrm{x},k(\mathrm{x}+\mathrm{e})) \end{array}\right]
\end{eqnarray}
The solution to~(\ref{cl_homog}) is denoted by \mbox{$\phi_t(x_0,e_0)=(\phi_t^x,\phi^e_t)^T$} (for $\phi_0(x_0,e_0)=(x_0,e_0)$), where $\phi_t^x$ represents the first $n$ coordinates of $\phi$ and $\phi_t^e$ represents the last $n$ coordinates. Under this notation, the execution rule $|\mathrm{e}|=c|\mathrm{x}|$ becomes $|\phi_t^e|=c|\phi_t^x|$. Since $f$ is homogeneous with respect to~(\ref{gen_dilation_vf}), the vector field~(\ref{cl_homog}) is homogeneous with respect to the following dilation vector field:
\begin{eqnarray}
%\label{gen_dilation_vff}
D = \sum_{i=1}^{n} r x_i \frac{\partial}{\partial x_i} + \sum_{i=1}^{n} r e_i \frac{\partial}{\partial e_i}
\end{eqnarray}
The homogeneous rays for this vector field are $\psi_s(x,e)= (\e^{rs}x,\e^{rs}e)$. %To clarify the argument, we define \mbox{$x_a= x_{t_i} (x_0)$}, the initial condition for the inter-sample behaviour at time $t_i$. 
%The flow $\phi$ depends on $(x_0,e_0)$, but since we are only interested in the case where $e_0=0$, we simply regard $\phi$ as a function of time and the initial condition of $x$. 
Since the closed loop~(\ref{cl_homog}) is homogeneous, $\phi_t$ satisfies~(\ref{Homog_traj}), i.e.: 
\begin{equation}
\phi_t (\e^{rs} x,\e^{rs} 0) = \e^{rs} \phi_{\e^{\rho(s)}t}(x,0), \qquad \forall x \in \R^n
\end{equation}
%This condition holds for all $t$, so \mbox{$\phi^x_{t_i} (\e^{rs} x) = \e^{rs} \phi^x_{\e^{\rho(s)}t_i} (x)$}, where $t_i$ represents the previous sampling instant. 
Hence if we consider now the initial condition of the dynamical system to be $(\e^{rs} x,\e^{rs} 0)$ the evolution of $|\phi^e_t|/|\phi^x_t|$ is:
\begin{eqnarray}
\label{time_rays}
\frac{\vert \phi^e_t     (\e^{rs} x,\e^{rs} 0)\vert}                  {\vert \phi^x_t(\e^{rs} x,\e^{rs} 0)\vert}  = 
%\frac{\vert \phi^x_{t_i} (\e^{rs} x)-\phi^x_t(\e^{rs} x)\vert}{\vert \phi^x_t(\e^{rs} x)\vert} \nonumber \\ & = & 
\frac{\vert \e^{rs}   \phi^e_{\e^{\rho(s)} t}(x,0)\vert}  {\vert \e^{rs}   \phi^x_{\e^{\rho(s)} t}(x,0)\vert}  \nonumber  = 
\frac{\vert \phi^e_{\e^{\rho(s)} t}(x,0)\vert}{\vert \phi^x_{\e^{\rho(s)} t}(x,0)\vert}
\end{eqnarray}
Therefore, the evolution of the ratio $|\phi^e_t|/|\phi^x_t|$ starting at $\e^{rs}(x,0)$ (the dilated point) is $\e^{\rho(s)}$ times faster than when starting at $(x,0)$: %Therefore, the inter-execution times will be $\e^{\rho(s)}$ shorter, as shown in equation~(\ref{scal_eq}): 
a dilation $\e^{rs}$ in the initial condition implies a dilation $\e^{\rho(s)}$ in the inter-execution times according to equation~(\ref{scal_eq_appendix}). %\qed
\end{proof}

\hspace{0.3cm} \textit{Theorem 5.1}: 
Consider the vector field $Z:M\rightarrow TM$, as defined in~(\ref{cl_gen}), and a vector field \mbox{$Y:N \rightarrow TN$} giving rise to flows $\mathrm{z}$ and $\mathrm{y}$, respectively. If there exist maps $\varphi:M \rightarrow N$ and $\eta_Y:N \rightarrow \R$ such that the following diagram commutes:

\begin{equation}
\label{comm_diag_app}
\begin{diagram}
\R  & \lTo^{\eta_Z}         & M & \rTo^Z         & TM\\
    & \luTo_{\eta_Y}        & \dTo^{\varphi}     &   & \dTo_{T\varphi}\\  
    &  & N & \rTo^Y & TN
\end{diagram}
\end{equation}
that is:
\begin{enumerate}
	\item $T\varphi Z = Y \circ \varphi$;
	\item $\eta_Y \circ \varphi = \eta_Z$,
\end{enumerate}
then the inter-execution times $\tau_z(z)$ implicitly defined by $\eta_Z \circ \mathrm{z}(t,z)= c$ (with $c>0$) and the inter-execution times $\tau_y(\varphi(z))$ implicitly defined by $ \eta_Y \circ \mathrm{y} (t,\varphi(z)) = c$ coincide:
\begin{equation}
\tau_z (z) = \tau_y (\varphi(z))
\end{equation}

\begin{proof}
Since $Z$ and $Y$ are $\varphi$-related, the corresponding flows satisfy~\cite{lee}:
$$\varphi \circ \mathrm{z}(t,z) = \mathrm{y} (t,\varphi(z))$$
and from commutativity of Diagram~\ref{comm_diag_app}, we can conclude:
\begin{equation}
\eta_Z \circ \mathrm{z}(t,z) = \eta_Y \circ \varphi \circ \mathrm{z}(t,z)  = \eta_Y \circ \mathrm{y} (t, \varphi(z))
\end{equation}
Hence, since both $\eta_Z \circ \mathrm{z}$ and $\eta_Y \circ \mathrm{y}$ are identical for any point $z \in M$, they generate the same sequence of times.
\end{proof}

\hspace{0.3cm} \textit{Corollary 5.2}: Let $\tau(x)$ be the inter-execution times for the system~(\ref{cl_gen}) implicitly defined by $|\mathrm{e}|=c|\mathrm{x}|$ for a point $x\in\R^n$, and let $\tilde{\tau}(x,w)$ be the inter-execution times for the system~(\ref{f_rel_vf}) under the same execution rule, for a point $(x,w)\in\R^{n+1}$. Then, $\tau$ and $\tilde \tau$ are related according to:%satisfies:
\begin{equation}
%\tau(x) = \lambda^{-l+1} \:\:\tilde{\tau} (x/\lambda,1/\lambda), \qquad \forall \lambda > 0
\tau(x) = \lambda^{l-1} \:\:\tilde{\tau} (\lambda x,\lambda), \qquad \forall \lambda > 0
\end{equation}
with $l$ as defined in~(\ref{defs_degrees}).

\begin{proof}
The system~(\ref{cl_gen}) is $\varphi$-related to~(\ref{f_rel_vf}), for $\varphi(x,e)=(x,1,e,0)$. If we define the output map $\eta_Y(x,w,e,e_w)=\frac{|e|}{|x|}$, then by Theorem~\ref{thm_phi_related} we can claim that the inter-execution times are identical:
$$ \tau(x) = \tilde \tau (x,1)$$
Moreover, since the auxiliary system~(\ref{f_rel_vf}) is homogeneous of degree $l-1$, times scale according to Theorem~\ref{ScaleHom}. Hence:
$$ \tau(x) = \tilde \tau (x,1) = \lambda^{l-1} \tilde \tau (\lambda x, \lambda), \qquad \forall \lambda >0$$

\end{proof}

\hspace{0.3cm} \textit{Proposition 5.3}: Consider the dynamical system~(\ref{cl_gen}), with $\hat f$ as defined in~(\ref{poly}). If $\Sigma$ is an invariant set for~(\ref{cl_gen}), then for any $\lambda>0$ the set $\tilde \Sigma_{\lambda} = \{(\lambda x,\lambda,\lambda e,0)\in \R^{2n+2}: (x,e)\in\Sigma\}$ is an invariant set for the extended system~(\ref{f_rel_vf}).	

\begin{proof}
Let the solution of~(\ref{cl_gen}) be denoted by $\mathrm{z}(t,x_0,e_0)=(\mathrm{x},\mathrm{e})^T$. Since $\Sigma$ is invariant, for any $(x_0,e_0) \in \Sigma$ the trajectory $\mathrm{z}$ stays in $\Sigma$, that is $\mathrm{z}(t,x_0,e_0) \in \Sigma$, for all $t \geq 0$. Now let the solution of the extended system~(\ref{f_rel_vf}) be denoted by $\mathrm{y}(t,x_0,w_0,e_0,e_{w_0})=(\tilde{ \mathrm{x}},\mathrm{w},\tilde{\mathrm{e}},\mathrm{e}_w)^T$. Using the fact that $\mathrm{y}$ is the flow of a homogeneous vector field, we can conclude:
\begin{eqnarray}
\lambda \mathrm{z}(t,x_0,e_0)&=&\lambda \left[ \begin{array}{ccc} \tilde{\mathrm{x}}(t,x_0,1,e_0,0)  \\ \tilde{\mathrm{e}}(t,x_0,1,e_0,0)  \end{array}\right] = 
\left[ \begin{array}{ccc} \tilde{\mathrm{x}}(\lambda^{l-1}t,\lambda x_0,\lambda,\lambda e_0,0) \\ \tilde{\mathrm{e}}(\lambda^{l-1}t,\lambda x_0,\lambda,\lambda e_0,0)  \end{array}\right]\notag\\
%\lambda \mathrm{x}(t,x_0,e_0)&=&\lambda\tilde{\mathrm{x}}(t,x_0,e_0,1,0) = \tilde{\mathrm{x}}(\lambda^{l-1}t,\lambda x_0,\lambda e_0,\lambda,0)\\
%\lambda \mathrm{e}(t,x_0,e_0)&=&\lambda\tilde{\mathrm{e}}(t,x_0,e_0,1,0) = \tilde{\mathrm{e}}(\lambda^{l-1}t,\lambda x_0,\lambda e_0,\lambda,0)
\end{eqnarray}
Moreover, we have $\mathrm{w}(t,\lambda x_0,\lambda , \lambda e_0,0)=\lambda$, for all $t\geq0$, as $\dot{\mathrm{w}}=0$. Hence if the initial condition of the auxiliary system~(\ref{f_rel_vf}) is $(\lambda x_0,\lambda ,\lambda e_0,0)$ for any $(x_0,e_0)\in\Sigma$, the trajectory $\mathrm{y}$ stays in $\tilde \Sigma_{\lambda}$ for $t\geq 0$. Thus every invariant set $\Sigma$ for the original system induces a family of invariant sets $\tilde \Sigma_{\lambda}$ for the auxiliary system defined as $\tilde \Sigma_{\lambda} = \{(\lambda x,\lambda,\lambda e,0) \in \R^{2n+2}: (x,e)\in\Sigma\}$, for any $\lambda>0$.

\end{proof}

\nocite{praly}	
\nocite{khalil2002ns}

\bibliographystyle{alpha}
\bibliography{self_triggered_control_nonlinear}

%\IEEEpeerreviewmaketitle

\end{document}